\documentclass[10pt]{article}
\usepackage{amssymb}
\usepackage{amsmath}
\usepackage[margin=1in]{geometry}

\catcode`\@=11
\@addtoreset{equation}{section}

\catcode`\@=12

\newtheorem{Theorem}{Theorem}[section]
\newtheorem{Definition}{Definition}[section]
\newtheorem{Lemma}{Lemma}[section]
\newtheorem{Corollary}{Corollary}[section]
\newtheorem{Remark}{Remark}[section]

\title{Lifespan of smooth solutions for timelike extremal surface equation in de Sitter spacetime}
\author{De-Xing Kong $\quad$and \quad Chang-Hua Wei\footnote{Corresponding
author: changhuawei1986@gmail.com.}\\
Department of Mathematics, Zhejiang University\\ Hangzhou 310027,
China}
\date{}
\begin{document}
\maketitle
\begin{abstract} In this paper, we study the generalized timelike extremal surface equation in the de Sitter spacetime, which plays an important role in both mathematics and physics. Under the assumption of small initial data with compact support, we investigate the lower bound of lifespan of smooth solutions by weighted energy estimates.

\vskip 3mm

\end{abstract}
\section{Introduction}
In this paper, we investigate the lifespan of smooth solutions for timelike extremal surface equation in the de Sitter spacetime under the assumption of small initial data. This kind of equation plays an important role in general relativity, the theory of black hole, particle physics, fluid mechanics and so on.
\subsection{Background}
The simplest family of black hole spacetimes with positive cosmological constant is the so-called Schwarzschild-de Sitter. If the cosmological constant $\Lambda>0$ is considered fixed, this is a 1-parameter family of the solution $(\mathcal{M},g)$ to the Einstein vacuum field equation
\begin{equation}
R_{\mu\nu}-\frac{1}{2}g_{\mu\nu}=-\Lambda g_{\mu\nu}
\end{equation}
with parameter $M$, called ``mass of black hole'', where $R_{\mu\nu}$ and $R$ are the Ricci curvature and scalar curvature of the manifold $\mathcal{M}$, respectively.

The line element has the form in local coordinates
\begin{equation}
ds^{2}=-(1-\frac{2M}{r}-\frac{1}{3}\Lambda r^{2})dt^{2}+(1-\frac{2M}{r}-\frac{1}{3}\Lambda r^{2})^{-1}dr^{2}+r^{2}d\sigma_{S^{2}},
\end{equation}
 where $d\sigma_{S^{2}}$ denotes the standard metric on the unit 2-sphere and $r$ is the standard distance in the Euclidean space.

In the present paper, we set $M=0$ to ignore the influence of the black hole. Then by the Lama$\hat{i}$tre-Robertson transformation \cite{M}, (1.2) reads
\begin{equation}
ds^{2}=-dt^{2}+e^{\frac{2t}{R}}(dx^{2}+dy^{2}+dz^{2}).
\end{equation}
The new coordinates $t,\,x,\,y,\,z$ can take all values from $-\infty$ to $\infty$ and $R$ is the ``radius'' of the universe. This is a special case of the line element of the Robertson-Walker space. For more details with respect to this spacetime, one can refer to \cite{H}.

The de Sitter line element in the higher dimensional analogue to the de Sitter space is
\begin{equation}
ds^{2}=-dt^{2}+e^{t}(\sum_{i=1}^{n}dx_{i}^{2}),
\end{equation}
where we have set $R=2$ for simplicity.

The timelike extremal surface equation corresponding to the de Sitter spacetime under the above coordinates (1.4) reads
\begin{equation}
\Box_{g}\varphi-\varphi_{t}=\frac{Q(\varphi,e^{t}Q(\varphi,\varphi))}{2(1+e^{t}Q(\varphi,\varphi))},
\end{equation}
where $\varphi=\varphi(t,x_{1},\cdots,x_{n})$ is the unknown function which corresponds to hypersurface, $Q(\varphi,\psi)$ is the null form of the de Sitter universe, which is given by
\begin{equation}
Q(\varphi,\psi)=-\varphi_{t}\psi_{t}+e^{-t}(\sum_{i=1}^{n}\varphi_{i}\psi_{i}).
\end{equation}
Here and hereafter, without confusion, we denote $\partial_{t}\varphi$ and $\partial_{x_{i}}\varphi$ by $\varphi_{t}$ and $\varphi_{i}$, respectively, $\Box_{g}$ denotes the covariant wave operator, which is given by
\begin{equation}
\Box_{g}=\frac{1}{\sqrt{|\det g|}}\partial_{\mu}(\sqrt{|\det g|}g^{\mu\nu}\partial_{\nu}),
\end{equation}
where the Einstein's summation convention has been used, namely, the same upper and lower index means summation. In the following, we will use this convention without a hint.

 Due to the important significances in both mathematics and physics, up to now, a lot of results on minimal surfaces in the Euclidean $\mathbb R^{n}$ and the Riemannian manifolds have been obtained. One can refer to two excellent books: Coding and Minicozzi \cite{Cod} and Osserman \cite{O}. Since the metric of Lorentz manifolds is not positive, a surface in these physical spacetimes may include the following four types: timelike, spacelike, lightlike and mixed types. For the global existence of extremal surface equation in Minkowski spacetime, one can refer to Lindblad \cite{Hans} and Brandle \cite{Br}. For the case in curved spaces, Gu \cite{Gu1} proves that in the isothermal coordinates, the equations of the motion of relativistic strings can be locally written as a form of harmonic map from Minkowski to the Lorentz manifold. He and Kong \cite{He} study the spherical solution of the relativistic membrane in Schwarzschild universe.

An important structure for the global existence of extremal surface equation in the Minkowski spacetime is the null condition, which is first studied by Klainerman \cite{K} and Christodoulou \cite{Ch}, respectively. For static curved spacetimes, Luk \cite{L2} investigates the global existence for nonlinear wave equations on slowly rotating Kerr spacetime satisfying the null condition, which is also suitable for the Schwarzschild spacetime. Other important models satisfying the null condition are the wave maps in curved spacetimes, which are interesting in geometry and physics, see \cite{B,Gu}. For more information on null conditions, one can refer to Alinhac \cite{A}. Obviously, the nonlinear terms of the extremal surface equation (1.5) satisfies the null condition, since $Q(\varphi,\psi)$ satisfies null condition by definition.

Recently, the wave equations in the background of de Sitter spacetime become the focus of interest for an increasing number of mathematicians. For linear wave equations, fundamental solutions to the Cauchy problem with or without source terms are obtained, see \cite{Y1,Y3}. For semilinear case, the global existence and blowup results for the Klein-Gordon equation have been arrived, see \cite{Y2,Y4}. To our knowledge, the present work is the first work on nonlinear wave equations in the background of de Sitter spacetime and we believe that it will play an important role in the study of the nonlinear stability of de Sitter universe.
\subsection{Main results}
Instead of considering (1.5), we shall consider Cauchy problems for the following wave equations
\begin{equation}
\Box_{g}\varphi-\varphi_{t}=Q(\varphi,\varphi)
\end{equation}
and
\begin{equation}
\Box_{g}\varphi-\varphi_{t}=\frac{Q(\varphi,e^{\alpha t}Q(\varphi,\varphi))}{2(1+e^{\alpha t}Q(\varphi,\varphi))},
\end{equation}
with the following initial data
\begin{equation}
t=0:\quad \varphi(0,x_{1},\cdots,x_{n})=\epsilon f(x_{1},\cdots,x_{n}),\; \varphi_{t}(0,x_{1},\cdots,x_{n})=\epsilon g(x_{1},\cdots,x_{n}),
\end{equation}
where $f,\,g\in C_{0}^{\infty}(\mathbb R^{n})$ and $\alpha\leq 1$ is a real number.
\begin{Remark}
When $\alpha=1$, (1.9) is nothing but (1.5).
\end{Remark}
Before we state our main results, we firstly give the definition of the lifespan of solution.
\begin{Definition}
The lifespan $T(\epsilon)$ is the supremum of $T>0$ such that the Cauchy problem (1.8), (1.10) (or (1.9), (1.10)) has a smooth solution on $[0,T]\times\mathbb R^{n}$.
\end{Definition}
\begin{Theorem}
There exists a positive constant $\epsilon_{0}$ such that, for any $\epsilon\in[0,\epsilon_{0}]$, the Cauchy problem (1.8), (1.10) has a unique global smooth solution on $[0,\infty)\times\mathbb R^{n}$.
\end{Theorem}
\begin{Theorem}
There exist positive constants $C_{0}$ and $\epsilon_{0}$, such that for $\epsilon\in[0,\epsilon_{0}]$, the lifespan of the Cauchy problem (1.9), (1.10) satisfies
\begin{equation}
T(\epsilon)=\left\{\begin{array}{ccc}&\infty,&\alpha<1,\\
&\frac{C_{0}}{\epsilon},&\alpha=1,
\end{array}
\right.
\end{equation}
where $C_{0}$ depends on $\epsilon_{0}$ and $\alpha$.
\end{Theorem}
\begin{Remark}
The results of Theorem 1.1 and 1.2 still hold if the dissipative term $\varphi_{t}$ does not appear in equations (1.8) and (1.9).
\end{Remark}
\subsection{Arrangement of the paper}
This paper is organized as follows. In Section 2, we investigate the basic equation for the motion of relativistic membrane in the de Sitter spacetime and derive an interesting nonlinear wave equation. In Section 3, we get a pointwise decay estimate for the linear wave equation with a dissipative term by the method of weighted energy estimates. In Section 4, we prove the global existence of the model equation (1.8). Section 5 is devoted to the lifespan of the Cauchy problem (1.9)-(1.10). Section 6 gives some discussions.
\section{Basic equation}
The de Sitter metric is given by (1.4), i.e.,
$$
ds^{2}=-dt^{2}+e^{t}(\sum_{i=1}^{n}dx_{i}^{2}).
$$
Consider the motion of a relativistic membrane in the de Sitter spacetime
$$
(t,x_{1},\cdots,x_{n})\rightarrow(t,x_{1},\cdots,x_{n},\varphi(t,x_{1},\cdots,x_{n})).
$$
In the coordinates $(t,x_{1},\cdots,x_{n})$, the induced metric of the submanifold $\mathcal{M}$ reads as
\begin{equation}
ds^{2}=(dt,dx_{1},\cdots,dx_{n})G(dt,dx_{1},\cdots,dx_{n})^{T},
\end{equation}
where
$$
G=\left(
\begin{array}{cccc}
g_{00}&g_{01}&\cdots&g_{0n}\\
g_{10}&g_{11}&\cdots&g_{1n}\\
\cdot\cdot&\cdot\cdot&\cdots&\cdot\cdot\\
g_{n0}&g_{n1}&\cdots&g_{nn}
\end{array}
\right),
$$
in which
$$
g_{00}=-1+e^{t}\varphi^{2}_{t},\quad g_{0i}=g_{i0}=e^{t}\varphi_{t}\varphi_{i}
$$
and
$$
g_{ii}=e^{t}+e^{t}\varphi^{2}_{i},\quad g_{ij}=e^{t}\varphi_{i}\varphi_{j}
$$
for $i,j=1,\cdots,n$.

We assume that the submanifold $\mathcal{M}$ is timelike, i.e.,
\begin{equation}
\Delta:=\det{G}=e^{nt}[-1-\sum_{i=1}^{n}\varphi_{i}^{2}+e^{t}\varphi_{t}^{2}]<0.
\end{equation}
This is equivalent to
$$
1+\sum_{i=1}^{n}\varphi_{i}^{2}-e^{t}\varphi_{t}^{2}>0.
$$
Thus the area element of $\mathcal{M}$ is
\begin{equation}
dA=\sqrt{-\Delta}dtdx_{1}\cdots dx_{n}.
\end{equation}
The submanifold $\mathcal{M}$ is called to be extremal if $\varphi=\varphi(t,x_{1},\cdots,x_{n})$ is a critical point of the area functional
\begin{equation}
I(\varphi)=\int\cdots\int\sqrt{-\Delta}dtdx_{1}\cdots dx_{n}.
\end{equation}
By direct calculations, the corresponding Euler-Lagrange equation reads as
\begin{equation}
\frac{\partial}{\partial t}\left(-e^{\frac{nt}{2}}\frac{e^{t}\varphi_{t}}{\sqrt{1+\sum_{i=1}^{n}(\varphi_{i}^{2})-e^{t}\varphi_{t}^{2}}}\right)+
\sum_{j=1}^{n}\frac{\partial}{\partial x_{j}}\left(e^{\frac{nt}{2}}\frac{\varphi_{j}}{\sqrt{1+\sum_{i=1}^{n}(\varphi_{i}^{2})-e^{t}\varphi_{t}^{2}}}\right),
\end{equation}
which is equivalent to
\begin{equation}
-\varphi_{tt}-\frac{n+2}{2}\varphi_{t}+e^{-t}\sum_{i=1}^{n}\varphi_{ii}+\frac{\varphi_{t}\partial_{t}(\sum_{i=1}^{n}\varphi_{i}^{2}
-e^{t}\varphi_{t}^{2})}{2(1+\sum_{i=1}^{n}\varphi_{i}^{2}-e^{t}\varphi_{t}^{2})}-
\sum_{j=1}^{n}e^{-t}\frac{\varphi_{j}\partial_{x_{j}}(\sum_{i=1}^{n}\varphi_{i}^{2}
-e^{t}\varphi_{t}^{2})}{2(1+\sum_{i=1}^{n}\varphi_{i}^{2}-e^{t}\varphi_{t}^{2})}=0.
\end{equation}
By direct calculatuons, the linear wave equation in de Sitter spacetime is
\begin{equation}
\Box_{g}\varphi=-\varphi_{tt}-\frac{n}{2}\varphi_{t}+e^{-t}\sum_{i=1}^{n}\varphi_{ii}.
\end{equation}
Denote
\begin{equation}
Q(\varphi,\psi)=-\varphi_{t}\psi_{t}+e^{-t}\sum_{i=1}^{n}\varphi_{i}\psi_{i},
\end{equation}
then by (2.7) and (2.8), (2.6) can be rewritten as
$$
\Box_{g}\varphi-\varphi_{t}=\frac{Q(\varphi,e^{t}Q(\varphi,\varphi))}{2(1+e^{t}Q(\varphi,\varphi))}.
$$
This is nothing but (1.5).

Before we state the structure enjoyed by (1.5), we generalize the definition of the null condition in the Minkowski spacetime to the de Sitter spacetime, which can be found in Alinhac \cite{A}.

Define $g^{\mu\nu}$ by
$$
g^{\mu\nu}g_{\mu\lambda}=\delta_{\lambda}^{\nu},
$$
where $\delta_{\lambda}^{\nu}$ is the Kronecker symbol.
\begin{Definition}
We say that a quadratic form
$$
A^{\mu\nu}\varphi_{\mu}\psi_{\nu}
$$
satisfies the null condition in general Lorentz manifold $(\mathcal{M},g)$, if the coefficients $A^{\mu\nu}$ satisfy
$$
A^{\mu\nu}\xi_{\mu}\xi_{\nu}=0,
$$
whenever $\xi$ is a null vector, namely, $ g(\xi,\xi)=g^{\mu\nu}\xi_{\mu}\xi_{\nu}=0$.
\end{Definition}
\begin{Lemma}
The nonlinear term $Q(\varphi,\psi)$satisfies the null condition in the sense of Definition (2.1).
\end{Lemma}
{\it Proof.}
It is easy to see that $Q(\varphi,\psi)=g^{\mu\nu}\varphi_{\mu}\psi_{v}$ satisfies the null condition. Thus, the lemma holds obviously by Definition 2.1.

Another important property for (2.5) is the linear degeneracy of its characteristics. In order to illustrate this phenomenon, we first recall the definition of linear degeneracy and genuine nonlinearity (see \cite{lax,lax1}).

Consider the following quasilinear hyperbolic systems
\begin{equation}
u_{t}+\sum_{k=1}^{n}A_{k}(u)u_{x_{k}}=B(u),
\end{equation}
where $u=(u_{1},\cdots,u_{n})^{T}$ is the unknown vector function, $A_{k}(u)=(a_{kij}(u))$ is an $n\times n$ matrix with suitably smooth elements $a_{kij}(u)$ $(i,j=1,\cdots,n)$, $B(u)=(B_{1}(u),\cdots,B_{n}(u))^{T}$ is a given smooth vector function, which denotes the source term. Define
\begin{equation}
A(u;\xi)=\sum_{k=1}^{n}A_{k}(u)\xi_{k},
\end{equation}
where $\xi=(\xi_{1},\cdots,\xi_{n})$ is any unit vector in the Euclidean space.

By hyperbolicity, for any given $u$ on the domain under consideration, $A(u;\xi)$ has $n$ real eigenvalues $\lambda_{1}(u;\xi),\cdots,\lambda_{n}(u;\xi)$ and a complete system of left (resp. right) eigenvectors. For $i=1,\cdots,n$,
let $l_{i}(u;\xi)=(l_{i1}(u;\xi),\cdots,l_{in}(u;\xi))$ (resp. $r_{i}(u;\xi)=(r_{i1}(u;\xi),\cdots,r_{in(u;\xi)})^{T}$) be a left (resp. right) eigenvector corresponding to $\lambda_{i}(u;\xi)$:
\begin{equation}
l_{i}(u;\xi)A(u;\xi)=\lambda_{i}(u;\xi)l_{i}(u;\xi)\quad(\text{resp.}\;A(u;\xi)r_{i}(u;\xi)=\lambda_{i}(u;\xi)r_{i}(u;\xi)).
\end{equation}
We have
\begin{equation}
\det|l_{ij}(u;\xi)|\neq 0\quad(\text{equivalently},\;\det|r_{ij}(u;\xi)|)=0
\end{equation}
Then
\begin{Definition}
$\lambda_{i}(u;\xi)$ $(i\in\{1,\cdots,n\})$ is said to be genuinely nonlinear, if for every state $u$ and any unit vector $\xi$, it holds that
\begin{equation}
\nabla\lambda_{i}(u;\xi)r_{i}(u;\xi)\neq0,
\end{equation}  $\lambda_{i}(u;\xi)$ is called to be linearly degenerate, if for every state $u$ and any unit vector $\xi$, it holds that
\begin{equation}
\nabla\lambda_{i}(u;\xi)r_{i}(u;\xi)\equiv0.
\end{equation}
The system (2.9) is genuinely nonlinear (resp. linearly degenerate), if all $\lambda_{i}$ $(i=1,\cdots,n)$ are genuinely nonlinear (resp. linearly degenerate).
\end{Definition}
Based on the above definition, we have
\begin{Lemma}
System (2.5) is linearly degenerate in the sense of P. D. Lax.
\end{Lemma}
{\it Proof.}
Set
\begin{equation}
\tau=-2e^{-\frac{t}{2}},\quad \tau\in[-2,0).
\end{equation}
We have
$$
\partial_{\tau}=\frac{dt}{d\tau}\partial_{t}=e^{\frac{t}{2}}\partial_{t},
$$
then
$$
\varphi_{\tau}=e^{\frac{t}{2}}\varphi_{t}.
$$
Thus, in the $(\tau,x_{1},\cdots,x_{n})$ coordinates, (2.5) can be rewritten  as
\begin{equation}
-\partial_{\tau}\left(\frac{\varphi_{\tau}}{\sqrt{1+\sum_{i=1}^{n}\varphi_{x_{i}}^{2}-\varphi_{\tau}^{2}}}\right)+\sum_{j=1}^{n}\partial_{x_{j}}\left(\frac{\varphi_{x_{j}}}{\sqrt{1+\sum_{i=1}^{n}\varphi_{x_{i}}^{2}-\varphi_{\tau}^{2}}}\right)=
-\frac{n+1}{\tau}\frac{\varphi_{\tau}}{\sqrt{1+\sum_{i=1}^{n}\varphi_{x_{i}}^{2}-\varphi_{\tau}^{2}}}.
\end{equation}
The principle term of (2.16) is nothing but the timelike extremal surface equation in the Minkowski spacetime $\mathbb R^{1+n}$, which is linearly degenerate obviously. One can refer to \cite{Kong}. Thus, the lemma is proved.

If the solution of (2.5) takes the following form
\begin{equation}
\varphi(t,x)=\varphi(t,x_{1},\cdots,x_{n}),
\end{equation}
where $x=\sum_{i=1}^{n}\xi_{i}x_{i}$ and $\xi=(\xi_{1},\cdots,\xi_{n})$ is the unit vector, then (2.5) can be reduced to
\begin{equation}
-\partial_{t}\left(\frac{e^{\frac{n+1}{2}t}e^{\frac{t}{2}}\varphi_{t}}{\sqrt{1+\varphi_{x}^{2}-e^{t}\varphi_{t}^{2}}}\right)+
\partial_{x}\left(\frac{e^{\frac{nt}{2}}\varphi_{x}}{\sqrt{1+\varphi_{x}^{2}-e^{t}\varphi_{t}^{2}}}\right)=0.
\end{equation}
Under the $(\tau, x)$ coordinate, (2.19) can be rewritten as
\begin{equation}
-\partial_{\tau}\left(\frac{\varphi_{\tau}}{\sqrt{1+\varphi_{x}^{2}-\varphi_{\tau}^{2}}}\right)+\partial_{x}\left(\frac{\varphi_{x}}{\sqrt{1+\varphi_{x}^{2}-\varphi_{\tau}^{2}}}\right)=
-\frac{n+1}{\tau}\frac{\varphi_{\tau}}{\sqrt{1+\varphi_{x}^{2}-\varphi_{\tau}^{2}}}.
\end{equation}
\begin{Remark}
In the $(\tau,x)$ coordinates, the principle term of (2.20) is nothing but the classical Born-Infeld equation \cite{In}.
\end{Remark}
\begin{Remark}
Different from the timelike extremal surface equation in the Minkowski spacetime $\mathbb R^{1+n}$, (2.17) has an extra dissipative term, since in the $(\tau,x_{1},\cdots,x_{n})$ coordinates, the coefficient of the source term $-\frac{n+1}{\tau}>0$, but the singularity appears as $\tau$ tends to zero.
\end{Remark}
\begin{Remark}
Equation (2.16) can be derived as the timelike extremal surface equation in the coordinates $(\tau,x_{1},\cdots,x_{n})$, where $\tau$ is defined by (2.15). In fact, in this coordinate frame, the metric (1.4) of de Sitter spacetime becomes
\begin{equation}
ds^{2}=\frac{4}{\tau^{2}}(-d\tau^{2}+\sum_{i=1}^{n}dx_{i}^{2}).
\end{equation}
\end{Remark}

\section{Pointwise decay for the linear wave equation}
In this section, we investigate the pointwise decay estimates of the following linear wave equation in de Sitter spacetime
\begin{equation}
\Box_{g}\varphi-\varphi_{t}=0.
\end{equation}
It will play a key role in the study of nonlinear cases. Before we state our main results of this section, we introduce the following notations
$$
\|u(x)\|_{L^{2}}=\left(\int_{\mathbb R^{n}}|u(x)|^{2}dx\right)^{\frac{1}{2}},\quad \|u(x)\|_{L^{\infty}}:=ess\,sup|u(x)|
$$
and
$$
\|u(x)\|_{H^{s}}=\left(\sum_{i=0}^{s}(\|D^{i}u(x)\|_{L^{2}})^{2}\right)^{\frac{1}{2}},
$$
where $s$ is an integer.

Define the energy momentum tensor corresponding to the equation (3.1) by
\begin{equation}
T_{\mu\nu}(\varphi)=\partial_{\mu}\varphi\partial_{\nu}\varphi-\frac{1}{2}g_{\mu\nu}|\nabla \varphi|^{2},
\end{equation}
where
$$
|\nabla \varphi|^{2}=g^{\kappa\lambda}\partial_{\kappa}\varphi\partial_{\lambda}\varphi=-\varphi_{t}^{2}+\sum_{i=1}^{n}e^{-t}\varphi_{i}^{2}.
$$
For a vector field $V=V^{\mu}\partial_{\mu}$, define the compatible currents
\begin{equation}
J_{\mu}^{V}(\varphi)=T_{\mu\nu}(\varphi)V^{\nu}
\end{equation}
and
\begin{equation}
K^{V}(\varphi)=\Pi_{\mu\nu}^{V}T^{\mu\nu}(\varphi),
\end{equation}
where $\Pi_{\mu\nu}^{V}$ is the deformation tensor defined by
\begin{equation}
\Pi_{\mu\nu}^{V}=\frac{1}{2}(\nabla_{\mu}V_{\nu}+\nabla_{\nu}V_{\mu}),
\end{equation}
in which $\nabla$ denotes the covariant derivative and
$$\nabla_{\mu}V_{\nu}=g(\nabla_{\mu}V,\partial_{\nu}).$$
For a constant $t$-slice, the induced volume form is defined by
\begin{equation}
dVol_{t}=e^{\frac{nt}{2}}dx_{1}\cdots dx_{n}.
\end{equation}
\begin{Remark}
In above notations, raising and lowering of indices in this paper is always done with respect to the metric $g$ of the de Sitter spacetime.
\end{Remark}

With above notations, by direct calculations, we have for $i,j=1,\cdots,n$ and $i\neq j$
\begin{equation}
T_{tt}(\varphi)=\varphi_{t}^{2}+\frac{1}{2}|\nabla\varphi|^{2}=\frac{1}{2}(\varphi_{t}^{2}+\sum_{i=1}^{n}e^{-t}\varphi_{i}^{2}),
\end{equation}
\begin{equation}
T_{ti}=\varphi_{t}\varphi_{i},\quad T_{ij}=\varphi_{i}\varphi_{j}
\end{equation}
and
\begin{equation}
T_{ii}=\varphi_{i}^{2}-\frac{1}{2}e^{t}|\nabla\varphi|^{2}=\frac{1}{2}(e^{t}\varphi_{t}^{2}+\varphi_{i}^{2}-\sum_{j\neq i}\varphi_{j}^{2}).
\end{equation}
The following lemma is easy and can be found in \cite{D1,D2,D3,D4,L1,L2}.
\begin{Lemma}
For the equation $\Box_{g}\varphi=f$, it holds that
\begin{equation}
\nabla^{\mu}T_{\mu\nu}=\Box_{g}\varphi\varphi_{\nu},\quad \nabla^{\mu}J_{\mu}^{V}(\varphi)=K^{V}(\varphi)+\Box_{g}\varphi\cdot V(\varphi).
\end{equation}
\end{Lemma}
{\it Proof.} By direct calculations, we have
\begin{align*}
\begin{split}
\nabla^{\mu}T_{\mu\nu}(\varphi)&=\nabla^{\mu}(\partial_{\mu}\varphi\partial_{\nu}\varphi-\frac{1}{2}g_{\mu\nu}\partial^{\lambda}\varphi\partial_{\lambda}\varphi)\\
&=\Box_{g}\varphi\partial_{\nu}\varphi+\partial_{\mu}\varphi\nabla^{\mu}\partial_{\nu}\varphi-g_{\mu\nu}\partial^{\lambda}\varphi\nabla^{\mu}\partial_{\lambda}\varphi\\
&=\Box_{g}\varphi\varphi_{\nu}
\end{split}
\end{align*}
and
\begin{align*}
\begin{split}
\nabla^{\mu}J_{\mu}^{V}(\varphi)&=\nabla^{\mu}(V^{\nu}T_{\mu\nu}(\varphi))\\
&=\nabla^{\mu}V^{\nu}T_{\mu\nu}(\varphi)+V^{\nu}\nabla^{\mu}T_{\mu\nu}(\varphi)\\
&=K^{V}(\varphi)+\Box_{g}\varphi V(\varphi).
\end{split}
\end{align*}
Thus, the lemma is proved.

The energy density $e(V,\upsilon)$ of the mapping $\varphi$ at time $t$ with respect to the past oriented timelike vector field $V$ is the nonnegative number
\begin{equation}
e(V,\upsilon)=J_{\alpha}^{V}\upsilon^{\alpha}=T_{\alpha\beta}(\varphi)V^{\beta}\upsilon^{\alpha}
\end{equation}
with $\upsilon^{\alpha}$ the components of the past oriented unit normal $\upsilon=-\partial_{t}$.

Taking the past oriented vector field $V$, by Lemma 3.1 and divergence theorem, we easily get the following lemma
\begin{Lemma}
The following energy identity holds in the domain $D=\{0\leq \tau\leq t\}$
\begin{equation}
\int_{\Sigma_{t}}J_{\alpha}^{V}\upsilon^{\alpha}dVol_{t}-\int_{\Sigma_{0}}J_{\alpha}^{V}\upsilon^{\alpha}dVol_{0}
=\int_{0}^{t}\int_{\Sigma_{\tau}}(K^{V}(\varphi)+\Box_{g}\varphi V(\varphi))dVol_{\tau}d\tau.
\end{equation}
\end{Lemma}

From now on, we take the vector field $V=-\partial_{t}$, then
\begin{equation}
\begin{split}
\Pi_{\mu\nu}^{V}&=\frac{1}{2}(\nabla_{\mu}V_{\nu}+\nabla_{\nu}V_{\mu})\\
&=\frac{1}{2}(g(\nabla_{\mu}(-\partial_{t}),\partial_{\nu})+g(\nabla_{\nu}(-\partial_{t}),\partial_{\mu}))\\
&=-\frac{1}{2}(g_{\nu\kappa}\Gamma_{\mu t}^{\kappa}+g_{\mu \kappa}\Gamma_{\nu t}^{\kappa}).
\end{split}
\end{equation}
So, for $i,j=1,\cdots,n$,
\begin{equation}
\Pi_{ii}^{V}=-g_{i\kappa}\Gamma_{i t}^{\kappa}=-g_{ii}\Gamma_{i t}^{i}=-\frac{1}{2}e^{t}.
\end{equation}
And for $i\neq j$,
\begin{equation}
\Pi_{ij}^{V}=0,\quad \Pi_{0i}^{V}=0\quad \text{and} \quad \Pi_{00}^{V}=0.
\end{equation}
Here $\Gamma_{ij}^{k}$ denotes the connection coefficients, which are given by
$$
\Gamma_{ij}^{k}=\frac{1}{2}g^{km}(\frac{\partial g_{im}}{\partial x_{j}}+\frac{\partial g_{jm}}{\partial x_{i}}-\frac{\partial g_{ij}}{\partial x_{m}}),
$$
where we have assumed that $x_{0}=t$.

By (3.4), (3.9), (3.14) and (3.15), we obtain
\begin{equation}
\begin{split}
K^{-\partial_{t}}(\varphi)&=\Pi_{\mu\nu}^{-\partial_{t}}T^{\mu\nu}(\varphi)\\
&=g^{\mu\mu}g^{\nu\nu}\Pi^{-\partial_{t}}_{\mu\nu}T_{\mu\nu}(\varphi)=g^{ii}g^{ii}\Pi_{ii}^{-\partial_{t}}T_{ii}(\varphi)\\
&=e^{-2t}\sum_{i=1}^{n}(-\frac{1}{2}e^{t})[\frac{1}{2}(e^{t}\varphi_{t}^{2}+\varphi_{i}^{2}-\sum_{j\neq i}\varphi_{j}^{2})]\\
&=\frac{1}{4}[(n-2)\sum_{i=1}^{n}e^{-t}\varphi_{i}^{2}-n\varphi_{t}^{2}].
\end{split}
\end{equation}
By (3.11),
\begin{equation}
e(V,\upsilon)=T_{\alpha\beta}(\varphi)V^{\alpha}\upsilon^{\beta}=T_{tt}(\varphi).
\end{equation}
Denote
$$
D=\{\partial_{1},\cdots,\partial_{n}\}\quad \text{and} \quad D^{I}=\partial_{1}^{I_{1}}\cdots\partial_{n}^{I_{n}},
$$
where $I=(I_{1},\cdots, I_{n})$ with $|I|=\sum_{j=1}^{n}|I_{j}|$.
For the constant $t$-slice, define
\begin{equation}
E^{|I|,I_{0}}_{0}(t)=\frac{1}{2}\int_{\Sigma_{t}}(\partial_{t}^{I_{0}}D^{I}\varphi)_{t}^{2}dVol_{t},\quad
 \quad
E^{|I|,I_{0}}_{1}(t)=\frac{1}{2}\int_{\Sigma_{t}}(\sum_{i=1}^{n}e^{-t}(\partial_{t}^{I_{0}}D^{I}\varphi)_{i}^{2})dVol_{t}
\end{equation}
and
\begin{equation}
E^{|I|,I_{0}}(t)=E^{|I|,I_{0}}_{0}(t)+E_{1}^{|I|,I_{0}}(t).
\end{equation}
Then, by above calculations and Lemma 3.2, the following zero-th order energy identity holds.
\begin{Lemma}
The energy identity (3.12) can be rewritten as
\begin{equation}
E^{0,0}(t)-E^{0,0}(0)=\int_{0}^{t}[-\frac{1}{2}(n+4)E^{0,0}_{0}(\tau)+\frac{1}{2}(n-2)E^{0,0}_{1}(\tau)]d\tau
\end{equation}
\end{Lemma}
{\it Proof.} By direct calculations, from (3.12), (3.16)-(3.18) and
$$
\int_{0}^{t}\int_{\Sigma_{\tau}}(-\varphi_{\tau}^{2})dVol_{\tau}d\tau=-2\int_{0}^{t}E^{0,0}_{0}(\tau)d\tau,
$$
we get (3.20) immediately.
\begin{Corollary}
It holds that
\begin{equation}
\frac{d}{dt}E^{0,0}(t)=-\frac{1}{2}(n+4)E^{0,0}_{0}(t)+\frac{1}{2}(n-2)E^{0,0}_{1}(t).
\end{equation}
\end{Corollary}

Based on the geometry of de Sitter spacetime with the metric given by (1.4), it is easy to see that the operator $D$ is a killing vector field, which means that $$\Pi_{\mu\nu}^{D}=0.$$ Thus, the structure of the equation (3.1) will not change if we take $D^{J}$ as a commutator, namely
\begin{equation}
\Box_{g}(D^{J}\varphi)-(D^{J}\varphi)_{t}=0.
\end{equation}
By (3.22) and Corollary 3.1, for $I_{0}=0$, we have
\begin{Corollary}
It holds that
\begin{equation}
\frac{d}{dt}E^{|J|,0}(t)=-\frac{1}{2}(n+4)E^{|J|,0}_{0}(t)+\frac{1}{2}(n-2)E^{|J|,0}_{1}(t),
\end{equation}
for arbitrary $J$.
\end{Corollary}
Define
\begin{equation}
f^{|I|,I_{0}}=E^{|I|,I_{0}}e^{\frac{2-n}{2}t},\quad f^{|I|,I_{0}}_{0}=E^{|I|,I_{0}}_{0}e^{\frac{2-n}{2}t}\quad\text{and}\quad f^{|I|,I_{0}}_{1}=E^{|I|,I_{0}}_{1}e^{\frac{2-n}{2}t}.
\end{equation}
Then we obtain
\begin{Lemma}
$f^{|I|,0}$ is uniformly bounded, provided that $E^{|I|,0}(0)$ is bounded for arbitrary $I$.
\end{Lemma}
{\it Proof.} By (3.24),
\begin{equation}
\begin{split}
\frac{d}{dt}f^{|I|,0}(t)&=\frac{d}{dt}E^{|I|,0}(t)e^{\frac{2-n}{2}t}+\frac{2-n}{2}E^{|I|,0}(t)e^{\frac{2-n}{2}t}\\
&=e^{\frac{2-n}{2}t}[-\frac{1}{2}(n+4)E^{|I|,0}_{0}(t)+\frac{1}{2}(n-2)E^{|I|,0}_{1}(t)+\frac{2-n}{2}(E^{|I|,0}_{0}(t)
+E^{|I|,0}_{1}(t))]\\
&=-(n-1)e^{\frac{2-n}{2}t}E^{|I|,0}_{0}(t)\leq 0.
\end{split}
\end{equation}
This proves Lemma 3.4.

For $I_{0}>0$, by (3.22), it holds that
\begin{Lemma}
\begin{equation}
\Box_{g}(\partial_{t}^{I_{0}}D^{J}\varphi)-(\partial_{t}^{I_{0}}D^{J}\varphi)_{t}=e^{-t}(\sum_{i=1}^{n}\sum_{M=0}^{I_{0}-1}C_{M}\partial_{t}^{M}\partial_{i}^{2}D^{J}\varphi),
\end{equation}
where $C_{M}$ $(M=0,\cdots,I_{0}-1)$ are constants depending on $M$.
\end{Lemma}
{\it Proof.} Denote
$$
D^{J}\varphi=v,
$$
it suffices to prove
\begin{equation}
\Box_{g}(\partial_{t}^{I_{0}}v)-(\partial_{t}^{I_{0}}v)_{t}=e^{-t}(\sum_{i=1}^{n}\sum_{M=0}^{I_{0}-1}C_{M}\partial_{t}^{M}\partial_{i}^{2}v).
\end{equation}
By (3.22), it holds that
\begin{equation}
\Box_{g}v-v_{t}=0.
\end{equation}
Since $\partial_{t}$ is not a killing vector field, it does not commutate with the operator $\Box_{g}$, by a direct calculation, it holds that
\begin{equation}
[\Box_{g},\partial_{t}]=[-\partial_{t}^{2}-\frac{n}{2}\partial_{t}+e^{-t}(\sum_{i=1}^{n}\partial_{i}^{2}),\partial_{t}]
=e^{-t}(\sum_{i=1}^{n}\partial_{i}^{2}).
\end{equation}
We prove this lemma by the method of induction.

When $I_{0}=1$, it holds that
\begin{equation}
\begin{split}
\Box_{g}(\partial_{t}v)&=[\Box_{g},\partial_{t}]v+\partial_{t}(\Box_{g}v)\\
&=e^{-t}(\sum_{i=1}^{n}\partial_{i}^{2}v)+(\partial_{t}v)_{t},
\end{split}
\end{equation}
i.e.,
\begin{equation}
\Box_{g}(\partial_{t}v)-(\partial_{t}v)_{t}=e^{-t}(\sum_{i=1}^{n}\partial_{i}^{2}v).
\end{equation}
Thus, the lemma holds for $I_{0}=1$.

Suppose the lemma holds for $I_{0}-1$, namely,
\begin{equation}
\Box_{g}(\partial_{t}^{I_{0}-1}v)-(\partial_{t}^{I_{0}-1}v)_{t}
=e^{-t}(\sum_{i=1}^{n}\sum_{M=0}^{I_{0}-2}C_{M}\partial_{t}^{M}\partial_{i}^{2}v),
\end{equation}
then
\begin{equation}
\begin{split}
\Box_{g}(\partial_{t}^{I_{0}}v)&=[\Box_{g},\partial_{t}](\partial_{t}^{I_{0}-1}v)+\partial_{t}(\Box_{g}(\partial_{t}^{I_{0}-1}v))\\
&=e^{-t}(\sum_{i=1}^{n}\partial_{i}^{2}\partial_{t}^{I_{0}-1}v)+\partial_{t}\left(
(\partial_{t}^{I_{0}-1}v)_{t}+e^{-t}(\sum_{i=1}^{n}\sum_{M=0}^{I_{0}-2}C_{M}\partial_{t}^{M}\partial_{i}^{2}v)\right)\\
&=e^{-t}(\sum_{i=1}^{n}\partial_{t}^{I_{0}-1}\partial_{i}^{2}v)+(\partial_{t}^{I_{0}}v)_{t}-e^{-t}
(\sum_{i=1}^{n}\sum_{M=0}^{I_{0}-2}C_{M}\partial_{t}^{M}\partial_{i}^{2}v)+e^{-t}(\sum_{i=1}^{n}\sum_{M=0}^{I_{0}-2}C_{M}
\partial_{t}^{M+1}\partial_{i}^{2}v)\\
&=(\partial_{t}^{I_{0}}v)_{t}+e^{-t}(\sum_{i=1}^{n}\sum_{M=0}^{I_{0}-1}C_{M}\partial_{t}^{M}\partial_{i}^{2}v),
\end{split}
\end{equation}
thus, the lemma holds for arbitrary $I_{0}$.

By Lemmas 3.2, 3.4, 3.5 and Corollary 3.2, for $I_{0}>0$, it holds that
\begin{Lemma}
$f^{|I|,I_{0}}(t)$ is uniformly bounded, and it holds that for arbitrary $I$,
\begin{equation}
f^{|I|,I_{0}}(t)\leq C_{I,I_{0}}(\sum_{k=0}^{I_{0}}\sum_{|l|+k\leq|I|+I_{0}}f^{|l|,k}(0)),
\end{equation}
where $C_{I,I_{0}}$ is a constant depending only on $I,\,I_{0}$.
\end{Lemma}
{\it Proof.} By Lemmas 3.2, 3.5 and Corollary 3.2, it is obvious that
\begin{equation}
\begin{split}
\frac{d}{dt}E^{|I|,I_{0}}(t)&=-\frac{1}{2}(n+4)E_{0}^{|I|,I_{0}}(t)+\frac{1}{2}(n-2)E_{1}^{|I|,I_{0}}(t)\\
&-
\int_{\Sigma_{t}}e^{-t}(\sum_{i=1}^{n}\sum_{M=0}^{I_{0}-1}C_{M}\partial_{t}^{M}\partial_{i}^{2}D^{|I|}\varphi)
\partial_{t}^{I_{0}+1}D^{|I|}\varphi dVol_{t}.
\end{split}
\end{equation}
As Lemma 3.4, it holds that
\begin{equation}
\begin{split}
\frac{d}{dt}f^{|I|,I_{0}}(t)&=\frac{d}{dt}E^{|I|,I_{0}}(t)e^{\frac{2-n}{2}t}+\frac{2-n}{2}E^{|I|,I_{0}}(t)e^{\frac{2-n}{2}t}\\
&\leq |e^{\frac{2-n}{2}t}\int_{\Sigma_{t}}e^{-t}(\sum_{i=1}^{n}\sum_{M=0}^{I_{0}-1}C_{M}\partial_{t}^{M}\partial_{i}^{2}D^{|I|}\varphi)
\partial_{t}^{I_{0}+1}D^{|I|}\varphi dVol_{t}|.
\end{split}
\end{equation}
Now, we prove the lemma by induction.

For $I_{0}=1$, by Lemma 3.4 and H\"{o}lder inequality, it holds that
\begin{equation}
\begin{split}
\frac{d}{dt}f^{|I|,1}(t)&\leq|e^{\frac{2-n}{2}t}\int_{\Sigma_{t}}e^{-t}(\sum_{i=1}^{n}\partial_{i}^{2}D^{|I|}\varphi)
\partial_{t}^{2}D^{|I|}\varphi dVol_{t}|\\
&\leq e^{-\frac{t}{2}}(f^{|I|+1,0}(t))^{\frac{1}{2}}(f^{|I|,1}(t))^{\frac{1}{2}}.
\end{split}
\end{equation}
Thus, by Lemma 3.4, it holds that
\begin{equation}
\left(f^{|I|,1}(t)\right)^{\frac{1}{2}}\leq \left(f^{|I|,1}(0)\right)^{\frac{1}{2}}+\int_{0}^{\infty}e^{-\frac{t}{2}}(f^{|I|+1,0}(0))^{\frac{1}{2}}dt,
\end{equation}
it implies that the lemma holds for $I_{0}=1$.

Suppose that the lemma holds for $N\leq I_{0}-1$, i.e.,
\begin{equation}
f^{|I|,N}(t)\leq C_{I,N}(\sum_{k=0}^{N}\sum_{|l|+k\leq|I|+N}f^{|l|,k}(0))\quad\text{for}\quad N\leq I_{0}-1,
\end{equation}
by (3.36) and H\"{o}lder inequality, we have
\begin{equation}
\frac{d}{dt}f^{|I|,I_{0}}(t)\leq \left(\sum_{M=0}^{I_{0}-1}C_{M}e^{-\frac{t}{2}}(f^{|I|+1,M}(t))^{\frac{1}{2}}\right)(f^{|I|,I_{0}}(t))^{\frac{1}{2}},
\end{equation}
thus, by (3.39) and (3.40), the lemma holds for arbitrary $I_{0}$.
\begin{Remark}
The quantities $f^{|I|,I_{0}}(0)$ and $E^{|I|,I_{0}}$ can be derived directly from the initial data and the equation.
\end{Remark}

Define
\begin{equation}
e^{|I|,I_{0}}_{0}(t)=\frac{1}{2}\|(\partial_{t}^{I_{0}}D^{I}\varphi)_{t}\|_{L^{2}}^{2}=\frac{1}{2}\int_{\mathbb R^{n}}(\partial_{t}^{I_{0}}D^{I}\varphi)_{t}^{2}dx_{1}\cdots dx_{n},
\end{equation}
\begin{equation}
e^{|I|,I_{0}}_{1}(t)=\frac{1}{2}\|\sum_{i=1}^{n}(\partial_{t}^{I_{0}}D^{I}\varphi)_{i}\|_{L^{2}}^{2}=\frac{1}{2}\int_{\mathbb R^{n}}\sum_{i=1}^{n}(\partial_{t}^{I_{0}}D^{I}\varphi)_{i}^{2}dx_{1}\cdots dx_{n}
\end{equation}
and
\begin{equation}
e^{|I|,I_{0}}(t)=e^{|I|,I_{0}}_{0}(t)+e^{|I|,I_{0}}_{1}(t).
\end{equation}
By (3.18), (3.41) and (3.42), we have
\begin{equation}
E^{|I|,I_{0}}_{0}(t)=\frac{1}{2}\int_{\mathbb R^{n}}(\partial_{t}^{I_{0}}D^{I}\varphi)_{t}^{2}e^{\frac{nt}{2}}dx_{1}\cdots dx_{n}=e^{\frac{nt}{2}}e^{|I|,I_{0}}_{0}(t)
\end{equation}
and
\begin{equation}
E^{|I|,I_{0}}_{1}(t)=\frac{1}{2}\int_{\mathbb R^{n}}\sum_{i=1}^{n}(\partial_{t}^{I_{0}}D^{I}\varphi)_{i}^{2}e^{-t}e^{\frac{nt}{2}}dx_{1}\cdots dx_{n}=e^{\frac{(n-2)t}{2}}e^{|I|,I_{0}}_{1}(t).
\end{equation}
Thus, we obtain easily
\begin{Lemma}
The following decay estimates hold
\begin{equation}
e^{|I|,I_{0}}_{0}(t)\leq e^{-t}f^{|I|,I_{0}}(t),\quad e^{|I|,I_{0}}_{1}(t)\leq f^{|I|,I_{0}}(t).
\end{equation}
\end{Lemma}

In what follows, we will use the following Sobolev embedding theorem on $\mathbb R^{n}$
\begin{Lemma}
If $u=u(x_{1},\cdots,x_{n})\in H^{s}$ for any $s>\frac{n}{2}$, then there exists a constant $C_{s}$ such that $$u\in L^{\infty}(\mathbb R^{n}),$$ and it holds that
\begin{equation}
\|u\|_{L^{\infty}}\leq C_{s}\|u\|_{H^{s}}.
\end{equation}
\end{Lemma}
By Lemmas (3.7) and (3.8), we easily obtain
\begin{Lemma}
For any $|I|\geq|J|+\lceil\frac{n}{2}+1\rceil$ and $i=1,\cdots,n$, it holds that
\begin{equation}
\|(\partial_{t}^{I_{0}}D^{J}\varphi)_{t}(t)\|_{L^{\infty}}\leq C_{I}e^{-\frac{t}{2}}(\sum_{|M|=0}^{|I|}f^{|M|,I_{0}}(t))^{\frac{1}{2}}\quad \text{for}\quad I_{0}\geq0
\end{equation}
and
\begin{equation}
\|(D^{J}\varphi)_{i}(t)\|_{L^{\infty}}\leq C_{I}(\sum_{|M|=0}^{|I|}f^{|M|,0}(t))^{\frac{1}{2}}.
\end{equation}
provided that $f^{|M|,I_{0}}(t)$ $(|M|=0,\cdots,|I|)$ is bounded. Here $\lceil a\rceil$ stands for the smallest integer larger than $a$.
\end{Lemma}
\begin{Remark}
By the discussions above, we observe that the dissipative term $\varphi_{t}$ does not affect the decay rate in this procedure.
\end{Remark}
\section{Global existence for Cauchy problem (1.8), (1.10)}
In this section, we shall consider the global existence of the following equation
$$
\Box_{g}\varphi=\varphi_{t}+Q(\varphi,\varphi),
$$
where $Q(\varphi,\varphi)$ is defined by (l.6) and satisfies Definition 2.1.

Since the equation can be reduced into a symmetric hyperbolic system, the local existence and uniqueness theorem holds, provided that the initial data belongs to the Sobolev space $H^{s}$ for $s>\frac{n}{2}+1$, which can be found in Majda \cite{Ma} and Alinhac \cite{A}.

The pointwise decay estimates derived in the last section will play a key role in the proof of the global existence and the lower bound of the lifespan for nonlinear wave equations. We will prove the main theorem by continity method and take the nonlinear terms as the disturbances. Before proving Theorem 1.1, we need the following lemmas, which state the structure enjoyed by null condition.
\begin{Lemma}
The null structure is conserved under $D$-derivatives, namely, the following holds
\begin{equation}
\Box_{g}(D^{I}\varphi)=(D^{I}\varphi)_{t}+\sum_{|I_{1}|+|I_{2}|=|I|}C_{I}Q(D^{I_{1}}\varphi,D^{I_{2}}\varphi),
\end{equation}
where $C_{I}$ is a constant depending on $I$.
\end{Lemma}
{\it Proof.}
Since the vector field $D$ is killing vector field corresponding to the operator $\Box_{g}$ and is commutable with $\partial_{t}$, it suffices to prove $$D^{I}Q(\varphi,\varphi)=\sum_{|I_{1}|+|I_{2}|=|I|}C_{I}Q(D^{I_{1}}\varphi,D^{I_{2}}\varphi).$$
In what follows, we prove it by induction.

When $|I|=1$, we have
$$
DQ(\varphi,\varphi)=D(-\varphi_{t}^{2}+\sum_{i=1}^{n}e^{-t}\varphi_{i}^{2})=-2\varphi_{t}(D\varphi)_{t}+\sum_{i=1}^{n}e^{-t}2\varphi_{i}(D\varphi)_{i}=2Q(D\varphi,\varphi).
$$
Thus, the lemma holds for $|I|=1$.

Suppose that the lemma holds for $|I|-1$, i.e., it holds that
$$
D^{I-1}Q(\varphi,\varphi)=\sum_{|\tilde{I}_{1}|+|\tilde{I}_{2}|=|I|-1}C_{I-1}Q(D^{\tilde{I}_{1}}\varphi,D^{\tilde{I}_{2}}\varphi)
,$$
then
\begin{equation}
\begin{split}
D^{I}Q(\varphi,\varphi)&=D(D^{I-1}Q(\varphi,\varphi))\\
&=D(\sum_{|\tilde{I}_{1}|+|\tilde{I}_{2}|=|I|-1}C_{I-1}Q(D^{\tilde{I}_{1}}\varphi,D^{\tilde{I}_{2}}\varphi))\\
&=D\left(\sum_{|\tilde{I}_{1}|+|\tilde{I}_{2}|=|I|-1}C_{I-1}(-(D^{\tilde{I}_{1}}\varphi)_{t}(D^{\tilde{I}_{2}}\varphi)_{t}+
\sum_{i=1}^{n}e^{-t}(D^{\tilde{I}_{1}}\varphi)_{i}(D^{\tilde{I}_{2}}\varphi)_{i})\right)\\
&=\sum_{|\tilde{I}_{1}|+|\tilde{I}_{2}|=|I|-1}C_{I-1}\left(-(DD^{\tilde{I}_{1}}\varphi)_{t}(D^{\tilde{I}_{2}}\varphi)_{t}+
\sum_{i=1}^{n}e^{-t}(DD^{\tilde{I}_{1}}\varphi)_{i}(D^{\tilde{I}_{2}}\varphi)_{i}\right)\\
&+\sum_{|\tilde{I}_{1}|+|\tilde{I}_{2}|=|I|-1}C_{I-1}\left(-(D^{\tilde{I}_{1}}\varphi)_{t}(DD^{\tilde{I}_{2}}\varphi)_{t}+
\sum_{i=1}^{n}e^{-t}(D^{\tilde{I}_{1}}\varphi)_{i}(DD^{\tilde{I}_{2}}\varphi)_{i}\right)\\
&=\sum_{|I_{1}|+|I_{2}|=|I|}C_{I}Q(D^{I_{1}}\varphi,D^{I_{2}}\varphi).
\end{split}
\end{equation}
Thus, the lemma holds.

For derivatives with respect to $t$, we have the following
\begin{Lemma}
It holds that
\begin{equation}\begin{split}
\Box_{g}(\partial_{t}^{I_{0}}D^{I}\varphi)-(\partial_{t}^{I_{0}}D^{I}\varphi)_{t}
&=e^{-t}(\sum_{i=1}^{n}\sum_{M=0}^{I_{0}-1}C_{M}\partial_{t}^{M}\partial_{i}^{2}D^{I}\varphi)\\
&+\sum_{\mbox{\tiny$\begin{array}{c}I_{01}+I_{02}=I_{0}\\|I_{1}|+|I_{2}|=|I|\end{array}$}}C_{I_{01},I_{02},I_{1},I_{2}}
Q(\partial_{t}^{I_{01}}D^{I_{1}}\varphi,\partial_{t}^{I_{02}}D^{I_{2}}\varphi)\\
&+e^{-t}(\sum_{i=1}^{n}\sum_{\mbox{\tiny$\begin{array}{c}
\tilde{I}_{01}+\tilde{I}_{02}\leq I_{0}-1\\|I_{1}|+|I_{2}|=|I|\end{array}$}}C_{\tilde{I}_{01},\tilde{I}_{02},I_{1},I_{2}}\partial_{t}^{\tilde{I_{01}}}\partial_{i}D^{I_{1}}
\varphi\partial_{t}^{\tilde{I_{02}}}\partial_{i}D^{I_{2}}\varphi),
\end{split}
\end{equation}
where $C_{M}$, $C_{I_{01},I_{02},I_{1},I_{2}}$ and $C_{\tilde{I}_{01},\tilde{I}_{02},I_{1},I_{2}}$ are constants.
\end{Lemma}
{\it Proof.} Denote
$$
D^{I_{1}}\varphi=v,\qquad D^{I_{2}}\varphi=w.
$$
By Lemmas 3.5 and 4.1, it suffices to prove
\begin{equation}
\partial_{t}^{I_{0}}Q(v,w)=\sum_{I_{01}+I_{02}=I_{0}}C_{I_{01},I_{02}}Q(\partial_{t}^{I_{01}}v,\partial_{t}^{I_{02}}w)
+e^{-t}(\sum_{i=1}^{n}\sum_{\tilde{I}_{01}+\tilde{I}_{02}\leq I_{0}-1}C_{\tilde{I}_{01},\tilde{I}_{02}}\partial_{t}^{\tilde{I}_{01}}\partial_{i}v\partial_{t}^{\tilde{I}_{02}}\partial_{i}w).
\end{equation}
As Lemmas 3.5 and 4.1, we prove (4.4) by induction.

When $I_{0}=1$, it holds that
\begin{equation}
\begin{split}
\partial_{t}Q(v,w)&=\partial_{t}\left(-\partial_{t}v\partial_{t}w+e^{-t}(\sum_{i=1}^{n}\partial_{i}v\partial_{i}w)\right)\\
&=-\partial_{t}^{2}v\partial_{t}w-\partial_{t}v\partial_{t}^{2}w+e^{-t}\left(\sum_{i=1}^{n}(\partial_{i}\partial_{t}v\partial_{i}w
+\partial_{i}v\partial_{i}\partial_{t}w)\right)-e^{-t}(\sum_{i=1}^{n}\partial_{i}v\partial_{i}w)\\
&=Q(\partial_{t}v,w)+Q(v,\partial_{t}w)-e^{-t}(\sum_{i=1}^{n}\partial_{i}v\partial_{i}w).
\end{split}
\end{equation}
Thus, the lemma holds for $I_{0}=1$.

Suppose the lemma holds for $I_{0}-1$, then
\begin{equation}
\begin{split}
\partial_{t}^{I_{0}}Q(v,w)&=\partial_{t}(\partial_{t}^{I_{0}-1}Q(v,w))\\
&=\partial_{t}\left(\sum_{I_{01}+I_{02}=I_{0}-1}C_{I_{01},I_{02}}Q(\partial_{t}^{I_{01}}v,\partial_{t}^{I_{02}}w)
+e^{-t}(\sum_{i=1}^{n}\sum_{\tilde{I}_{01}+\tilde{I}_{02}\leq I_{0}-2}\partial_{t}^{\tilde{I}_{01}}\partial_{i}v\partial_{t}^{\tilde{I}_{02}}\partial_{i}w)\right)\\
&=\sum_{I_{01}+I_{02}=I_{0}-1}C_{I_{01},I_{02}}\left(Q(\partial_{t}^{I_{01}+1}v,\partial_{t}^{I_{02}}w)+
Q(\partial_{t}^{I_{01}}v,\partial_{t}^{I_{02}+1}w)\right)\\
&-e^{-t}(\sum_{i=1}^{n}\partial_{t}^{I_{01}}\partial_{i}v\partial_{t}^{I_{02}}\partial_{i}w)
-e^{-t}(\sum_{i=1}^{n}\sum_{\tilde{I}_{01}+\tilde{I}_{02}\leq I_{0}-2}\partial_{t}^{\tilde{I}_{01}}\partial_{i}v\partial_{t}^{\tilde{I}_{02}}\partial_{i}w)\\
&+e^{-t}\left(\sum_{i=1}^{n}\sum_{\tilde{I}_{01}+\tilde{I}_{02}\leq I_{0}-2}\partial_{t}^{\tilde{I}_{01}+1}\partial_{i}v\partial_{t}^{\tilde{I}_{02}}\partial_{i}w
+\sum_{i=1}^{n}\sum_{\tilde{I}_{01}+\tilde{I}_{02}\leq I_{0}-2}\partial_{t}^{\tilde{I}_{01}}\partial_{i}v\partial_{t}^{\tilde{I}_{02}+1}\partial_{i}w\right)\\
&=\sum_{I_{01}+I_{02}=I_{0}}C_{I_{01},I_{02}}Q(\partial_{t}^{I_{01}}v,\partial_{t}^{I_{02}}w)
+e^{-t}(\sum_{i=1}^{n}\sum_{\tilde{I}_{01}+\tilde{I}_{02}\leq I_{0}-1}\partial_{t}^{\tilde{I}_{01}}\partial_{i}v\partial_{t}^{\tilde{I}_{02}}\partial_{i}w)
\end{split}
\end{equation}
Thus, the lemma holds for arbitrary $I_{0}$ and $I$.

Now we prove Theorem 1.1.\\
{\it Proof of Theorem 1.1.} We prove the theorem by the following three steps.

Step 1: Define
\begin{equation}
F(t)=\sum_{I_{0}+|I|\leq N}f^{|I|,I_{0}}(t) \quad 0\leq t< T,
\end{equation}
for $N\geq n+4$, where $f^{|I|,I_{0}}(t)$ is defined by (3.24) whenever $\varphi\in C^{\infty}([0,T)\times\mathbb R^{n})$ solves (1.8), (1.10) on $[0,T)\times\mathbb R^{n}$ for some $T>0$. By (1.10), there exists a positive constant $C_{0}$ depends only on the initial data $f$, $g$ and their derivatives such that
\begin{equation}
F(0)\leq C_{0}\epsilon.
\end{equation}
Then by the pointwise decay estimates of last section, when $N\geq max\{|J_{1}|+|J_{2}|-1,|J_{3}|-1\}+\lceil\frac{n}{2}+1\rceil$, it holds that
\begin{equation}
\|\partial_{t}^{J_{1}}D^{J_{2}}\varphi\|_{L^{\infty}}\leq Ce^{-\frac{t}{2}}F^{\frac{1}{2}}(t), \quad \text{when}\quad|J_{1}|\geq1
\end{equation}
and
\begin{equation}
\|D^{J_{3}}\varphi\|_{L^{\infty}}\leq CF^{\frac{1}{2}}(t),
\end{equation}
where $C$ is a constant coming from the Sobolev embedding theorem.

Step 2: Energy estimates\\
Suppose $|J|+M\leq N$, by Corollary 3.2, Lemmas 3.9, 4.2 and (3.36), we have
\begin{equation}
\begin{split}
&\frac{d}{dt}f^{|J|,M}(t)\\
&\leq
|e^{\frac{2-n}{2}t}\int_{\Sigma_{t}}e^{-t}(\sum_{i=1}^{n}\sum_{K=0}^{M-1}C_{K}\partial_{t}^{K}\partial_{i}^{2}D^{|J|}\varphi)
\partial_{t}^{M+1}D^{|J|}\varphi dVol_{t}|\\
&+|e^{\frac{2-n}{2}t}\int_{\Sigma_{t}}\sum_{\mbox{\tiny$\begin{array}{c}I_{01}+I_{02}=M\\|J_{1}|+|J_{2}|=|J|\end{array}$}}C_{M_{01},M_{02},J_{1},J_{2}}
Q(\partial_{t}^{M_{01}}D^{J_{1}}\varphi,\partial_{t}^{M_{02}}D^{J_{2}}\varphi)\partial_{t}^{M+1}D^{|J|}\varphi dVol_{t}|\\
&+|e^{\frac{2-n}{2}t}\int_{\Sigma_{t}}e^{-t}(\sum_{i=1}^{n}\sum_{\mbox{\tiny$\begin{array}{c}
\tilde{M}_{01}+\tilde{M}_{02}\leq M-1\\|J_{1}|+|J_{2}|=|J|\end{array}$}}C_{\tilde{M}_{01},\tilde{M}_{02},J_{1},J_{2}}\partial_{t}^{\tilde{M}_{01}}\partial_{i}D^{J_{1}}
\varphi\partial_{t}^{\tilde{M_{02}}}\partial_{i}D^{J_{2}}\varphi)\partial_{t}^{M+1}D^{|J|}\varphi dVol_{t}|\\
&\leq e^{-\frac{t}{2}}\left(\sum_{K=0}^{M-1}|C_{K}|(f^{|J|+1,K}(t))^{\frac{1}{2}}(f^{|J|,M}(t))^{\frac{1}{2}}\right)\\
&+\sum_{\mbox{\tiny$\begin{array}{c}M_{01}+M_{02}=M\\|J_{1}|+|J_{2}|=|J|\end{array}$}}C_{M_{01},M_{02},J_{1},J_{2}}
(\|\partial_{t}^{M_{01}+1}D^{J_{1}}\varphi\|_{L^{\infty}}+e^{-\frac{t}{2}}\|\partial_{t}^{M_{01}}\partial_{i}D^{J_{1}}\varphi\|_{L^{\infty}})\times\\
&\quad\left(f^{|J_{2}|,M_{02}}(t)\right)^{\frac{1}{2}}
\left(f^{|J|,M}(t)\right)^{\frac{1}{2}}\\
&+\sum_{i=1}^{n}\sum_{\mbox{\tiny$\begin{array}{c}
\tilde{M}_{01}+\tilde{M}_{02}\leq M-1\\|J_{1}|+|J_{2}|=|J|\end{array}$}}C_{\tilde{M}_{01},\tilde{M}_{02},J_{1},J_{2}}
e^{-\frac{t}{2}}\|\partial_{t}^{\tilde{M}_{01}}\partial_{i}D^{J_{1}}\varphi\|_{L^{\infty}}\left(f_{1}^{|J_{2}|,\tilde{M}_{02}}(t)\right)^{\frac{1}{2}}
\left(f^{|J|,M}(t)\right)^{\frac{1}{2}}
\end{split}
\end{equation}
where we have assumed without loss of generality
$$M_{01}+|J_{1}|\leq M_{02}+|J_{2}|\quad\text{and}\quad \tilde{M}_{01}+|J_{1}|\leq\tilde{M}_{02}+|J_{2}|$$
 Summing $|J|$ and $M$, which satisfy $|J|+M\leq N$, we easily get that
\begin{equation}
\frac{d}{dt}F(t)\leq C_{N}(e^{-\frac{t}{2}}+\|\partial_{t}^{M_{01}+1}D^{J_{1}}\varphi\|_{L^{\infty}}
+e^{-\frac{t}{2}}\|\partial_{t}^{M_{01}}\partial_{i}D^{J_{1}}\varphi\|_{L^{\infty}}
+e^{-\frac{t}{2}}\|\partial_{t}^{\tilde{M}_{01}}\partial_{i}D^{J_{1}}\varphi\|_{L^{\infty}})F(t),
\end{equation}
where
$$
M_{01}+|J_{1}|\leq\frac{N}{2}\quad \text{and}\quad \tilde{M}_{01}+|J_{1}|\leq\frac{N-1}{2}.
$$
By Lemma 3.9, if $M_{01}+|J_{1}|+1+\frac{n}{2}+1\leq\frac{N}{2}+\frac{n}{2}+2\leq N$, i.e., $N\geq n+4$, it holds that
\begin{equation}
\frac{d}{dt}F(t)\leq C_{N}(e^{-\frac{t}{2}}+e^{-\frac{t}{2}}F^{\frac{1}{2}}(t))F(t).
\end{equation}

Step 3: Boot-strap \\
Set
$$
E=\{t\in[0,T):F(s)\leq A\epsilon\quad\text{for}\,\text{all}\quad0\leq s\leq t\}.
$$
By (4.8), $E$ is not empty. Since $F(t)$ is continuous in $t$, $E$ is relatively closed in $[0,T)$. Thus, it suffices to prove that $E$ is relatively open such that the following holds.

For any $T$, set
$$
E=[0,T).
$$

 In order to prove $E$ is open, we fix $t_{0}\in E$ with $t_{0}<T$. Since $F(t)$ is continuous, there exists $t_{1}>t_{0}$ such that
\begin{equation}
F(t)\leq 2A\epsilon\quad\text{for}\quad 0\leq t\leq t_{1}.
\end{equation}
We shall prove
\begin{equation}
F(t)\leq A\epsilon\quad\text{for}\quad 0\leq t\leq t_{1},
\end{equation}
provided $\epsilon$ is sufficiently small. By (4.13), in the domain $[0,t_{1}]$, we have
\begin{equation}
\frac{d}{dt}F(t)\leq C_{N}(e^{-\frac{t}{2}}+e^{-\frac{t}{2}}(2A\epsilon)^{\frac{1}{2}})F(t).
\end{equation}
Thus, we obtain
\begin{equation}
F(t)\leq F(0)e^{\int_{0}^{t}C_{N}(e^{-\frac{t}{2}}+e^{-\frac{t}{2}}(2A\epsilon)^{\frac{1}{2}})dt}\leq C_{0}\epsilon
e^{\int_{0}^{\infty}C_{N}(e^{-\frac{t}{2}}+e^{-\frac{t}{2}}(2A\epsilon)^{\frac{1}{2}})dt}.
\end{equation}
If $\epsilon$ is sufficiently small, then the following holds obviously
\begin{equation}
e^{\int_{0}^{\infty}C_{N}(e^{-\frac{t}{2}}+e^{-\frac{t}{2}}(2A\epsilon)^{\frac{1}{2}})dt}\leq \frac{A}{C_{0}},
\end{equation}
provided $A$ is sufficiently large. Thus, we can get that the solution of the Cauchy problem (1.8), (1.10) exists globally by the standard continuity method. This completes the proof.
\begin{Remark}
From the above procedure, we can see that the null condition $Q(\varphi,\varphi)$ plays a key role, especially the coefficients $e^{-t}$ before $\varphi_{i}^{2}$. In fact, $e^{-\frac{t}{2}-\delta}$ is enough for the global existence result of (1.8), provided that $\delta>0$. In next section, we will clarify the influence of this term.
\end{Remark}
\section{Lifespan for Cauchy problem (1.9)-(1.10)}
In this section, we consider the more complicated and representative case (1.9), which generalizes the timelike extremal surface equation and is interesting in the fields of both mathematics and physics. Since the nonlinearity is higher than the components of (1.8), we must generalize the corresponding energy estimates. As before, we study the structure enjoyed by the nonlinear term after differentiated several times by $D$.

Define
\begin{equation}
S(I):=\left\{\begin{array}{cc}
(I_{1},I_{2},I_{3},I_{4_{m}},I_{5_{m}}):|I_{1}|+|I_{2}|+|I_{3}|+|I_{4_{m}}|+|I_{5_{m}}|\leq|I|\\
\max\{|I_{2}|,|I_{3}|\}\leq|I|-1,\;\text{for}\;0\leq j\leq |I|,\,0\leq m\leq j
\end{array}
\right\}.
\end{equation}
\begin{Lemma}
Differentiate the nonlinear term $\frac{Q(\varphi,e^{\alpha t}Q(\varphi,\varphi))}{2(1+e^{\alpha t}Q(\varphi,\varphi))}$ for $|I|$ times by $D$, we have
\begin{equation}
\begin{split}
&D^{I}\frac{Q(\varphi,e^{\alpha t}Q(\varphi,\varphi))}{2(1+e^{\alpha t}Q(\varphi,\varphi))}=\frac{Q(\varphi,e^{\alpha t}Q(\varphi,D^{I}\varphi))}{(1+e^{\alpha t}
Q(\varphi,\varphi))}\\
&\quad+\sum_{j=0}^{|I|}\sum_{S(I)}\prod_{m=0}^{j}G(e^{ \alpha t}Q(\varphi,\varphi))Q(D^{I_{1}}\varphi,e^{ \alpha t}Q(D^{I_{2}}\varphi,D^{I_{3}}\varphi))
e^{j\alpha t}Q(D^{I_{4_{m}}}\varphi,D^{I_{5_{m}}}\varphi),
\end{split}
\end{equation}
where $G(e^{\alpha t}Q(\varphi,\varphi))$ is a smooth function depending on $e^{\alpha t}Q(\varphi,\varphi)$ and when $m=0$, the terms containing the index $m$ in the product terms do not appear.
\end{Lemma}
{\it Proof.} The proof is by induction on $I$ and using Lemma 4.1 repeatedly.

When $|I|=1$, we have
\begin{equation}\begin{split}
&D\frac{Q(\varphi,e^{\alpha t}Q(\varphi,\varphi))}{2(1+e^{\alpha t}Q(\varphi,\varphi))}\\
& =\frac{Q(\varphi,e^{\alpha t}Q(\varphi,D\varphi))}{(1+e^{\alpha t}Q(\varphi,\varphi))}+
\frac{Q(D\varphi,e^{\alpha t}Q(\varphi,\varphi))}{2(1+e^{\alpha t}Q(\varphi,\varphi))}
-\frac{Q(\varphi,e^{\alpha t}Q(\varphi,\varphi))e^{\alpha t}Q(\varphi,D\varphi)}{2(1+e^{\alpha t}Q(\varphi,\varphi))^{2}}\\
&=\frac{Q(\varphi,e^{\alpha t}Q(\varphi,D\varphi))}{(1+e^{\alpha t}Q(\varphi,\varphi))}+\sum_{j=0}^{1}\sum_{S(1)}G(e^{\alpha t}Q(\varphi,\varphi))Q(D^{I_{1}}\varphi,e^{\alpha t}Q(D^{I_{2}}\varphi,D^{I_{3}}\varphi))
e^{j\alpha t}Q(D^{I_{4_{j}}}\varphi,D^{I_{5_{j}}}\varphi).
\end{split}
\end{equation}
Thus, the lemma holds for $|I|=1$.

Suppose that the lemma holds for $|J|=|I|-1$, then for $|I|$
\begin{equation}
\begin{split}
&D^{I}\frac{Q(\varphi,e^{\alpha t}Q(\varphi,\varphi))}{2(1+e^{\alpha t}Q(\varphi,\varphi))}=D\left(D^{I-1}\frac{Q(\varphi,e^{\alpha t}Q(\varphi,\varphi))}{2(1+e^{\alpha t}Q(\varphi,\varphi))}\right)=D\frac{Q(\varphi,e^{\alpha t}Q(\varphi,D^{J}\varphi))}{(1+e^{\alpha t}
Q(\varphi,\varphi))}\\
&+D\left(\sum_{j=0}^{|J|}\sum_{S(J)}\prod_{m=0}^{j}G(e^{\alpha t}Q(\varphi,\varphi))Q(D^{J_{1}}\varphi,e^{\alpha t}Q(D^{J_{2}}\varphi,D^{J_{3}}\varphi))
e^{j\alpha t}Q(D^{J_{4_{m}}}\varphi,D^{J_{5_{m}}}\varphi)\right)\\
&=\frac{Q(\varphi,e^{\alpha t}Q(\varphi,D^{I}\varphi))}{1+e^{\alpha t}Q(\varphi,\varphi)}+
\frac{Q(D\varphi,e^{\alpha t}Q(\varphi,D^{J}\varphi))}{1+e^{\alpha t}Q(\varphi,\varphi)}\\
&+
\frac{Q(\varphi,e^{\alpha t}Q(D\varphi,D^{J}\varphi))}{1+e^{\alpha t}Q(\varphi,\varphi)}
-
\frac{2Q(\varphi,e^{\alpha t}Q(\varphi,D^{J}\varphi))e^{\alpha t}Q(\varphi,D\varphi)}{(1+e^{\alpha t}Q(\varphi,\varphi))^{2}}\\
&+\sum_{j=0}^{|J|}\sum_{S(J)}D\left(\prod_{m=0}^{j}G(e^{\alpha t}Q(\varphi,\varphi))Q(D^{J_{1}}\varphi,e^{\alpha t}Q(D^{J_{2}}\varphi,D^{J_{3}}\varphi))
e^{j\alpha t}Q(D^{J_{4_{m}}}\varphi,D^{J_{5_{m}}}\varphi)\right)\\
&=\frac{Q(\varphi,e^{\alpha t}Q(\varphi,D^{I}\varphi))}{1+e^{\alpha t}Q(\varphi,\varphi)}\\
&+\sum_{j=0}^{|I|}\sum_{S(I)}\prod_{m=0}^{j}G(e^{\alpha t}Q(\varphi,\varphi))Q(D^{I_{1}}\varphi,e^{\alpha t}Q(D^{I_{2}}\varphi,D^{I_{3}}\varphi))
e^{j\alpha t}Q(D^{I_{4_{m}}}\varphi,D^{I_{5_{m}}}\varphi).
\end{split}
\end{equation}
The last equality comes from the product role. Thus, the lemma is proved.

For derivatives with respect to $t$, the following lemmas play key roles.
\begin{Lemma}
It holds that
\begin{equation}
\begin{split}
\partial_{t}^{J}Q(u,e^{\alpha t}Q(v,w))&=\sum_{J_{1}+J_{2}+J_{3}\leq J}C(J_{1},J_{2},J_{3})Q(\partial_{t}^{J_{1}}u,e^{\alpha t}Q(\partial_{t}^{J_{2}}v,\partial_{t}^{J_{3}}w))\\
&+\sum_{i=1}^{n}\sum_{J_{4}+J_{5}+J_{6}\leq J-1}C(J_{4},J_{5},J_{6})Q(\partial_{t}^{J_{4}}u,e^{(\alpha-1)t}\partial_{t}^{J_{5}}v_{i}\partial_{t}^{J_{6}}w_{i})\\
&+\sum_{i=1}^{n}\sum_{J_{7}+J_{8}+J_{9}\leq J-1}C(J_{7},J_{8},J_{9})e^{(\alpha-1)t}\partial_{t}^{J_{7}}u_{i}\partial_{i}\left(Q(\partial_{t}^{J_{8}}v,\partial_{t}^{J_{9}}w)\right)\\
&+\sum_{i,j=1}^{n}\sum_{J_{10}+J_{11}+J_{12}\leq J-2}C(J_{10},J_{11},J_{12})e^{(\alpha-2)t}\partial_{t}^{J_{10}}u_{i}\partial_{t}^{J_{11}}v_{j}\partial_{t}^{J_{12}}w_{j},\\
\end{split}
\end{equation}
where $J,\,J_{i}\,(i=1,\cdots,12)$ are non-negative integers and $C(J)$ denotes a constant depending on $J$.
\end{Lemma}
{\it Proof.} For simplicity, we neglect the constants in the proof. By (4.4), it holds that
\begin{equation}
\begin{split}
&\partial_{t}^{J}Q(u,e^{\alpha t}Q(v,w))\\
&=\sum_{a_{1}+a_{2}=J}Q\left(\partial_{t}^{a_{1}}u,\partial_{t}^{a_{2}}(e^{\alpha t}Q(v,w))\right)+
\sum_{i=1}^{n}\sum_{a_{3}+a_{4}\leq J-1}e^{-t}\partial_{t}^{a_{3}}u_{i}\partial_{t}^{a_{4}}\partial_{i}(e^{\alpha t}Q(v,w))\\
&=\sum_{a_{1}+a_{2}=J}Q\left(\partial_{t}^{a_{1}}u,\sum_{a_{5}+a_{6}=a_{2}}\partial_{t}^{a_{5}}(e^{\alpha t})\partial_{t}^{a_{6}}Q(v,w)\right)\\
&+\sum_{i=1}^{n}\sum_{a_{3}+a_{4}\leq J-1}e^{-t}\partial_{t}^{a_{3}}u_{i}\sum_{a_{7}+a_{8}=a_{4}}\partial_{t}^{a_{7}}(e^{\alpha t})\partial_{t}^{a_{8}}\partial_{i}Q(v,w)\\
&=\sum_{a_{1}+a_{2}=J}Q\left(\partial_{t}^{a_{1}}u,\sum_{a_{5}+a_{6}=a_{2}}\sum_{a_{9}+a_{10}=a_{6}}\partial_{t}^{a_{5}}(e^{\alpha t})Q(\partial_{t}^{a_{9}}v,\partial_{t}^{a_{10}}w)\right)\\
&+\sum_{a_{1}+a_{2}=J}Q\left(\partial_{t}^{a_{1}}u,\sum_{i=1}^{n}\sum_{a_{11}+a_{12}\leq a_{6}-1}\sum_{a_{5}+a_{6}=a_{2}}\partial_{t}^{a_{5}}(e^{\alpha t})e^{-t}\partial_{t}^{a_{11}}v_{i}\partial_{t}^{a_{12}}w_{i}\right)\\
&+\sum_{i=1}^{n}\sum_{a_{3}+a_{4}\leq J-1}e^{-t}\partial_{t}^{a_{3}}u_{i}\sum_{a_{7}+a_{8}=a_{4}}\partial_{t}^{a_{7}}(e^{\alpha t})\partial_{i}\left(\sum_{a_{13}+a_{14}=a_{8}}Q(\partial_{t}^{a_{13}}v,\partial_{t}a^{14}w)\right)\\
&+\sum_{i=1}^{n}\sum_{a_{3}+a_{4}\leq J-1}e^{-t}\partial_{t}^{a_{3}}u_{i}\sum_{a_{7}+a_{8}=a_{4}}\partial_{t}^{a_{7}}(e^{\alpha t})\partial_{i}\left(\sum_{j=1}^{n}\sum_{a_{15}+a_{16}\leq a_{8}-1}e^{-t}\partial_{t}^{a_{15}}v_{j}\partial_{t}^{a_{16}}w_{j}\right).
\end{split}
\end{equation}
Rearrange the indices $a_{i}\,(i=1,\cdots,16)$ of (5.6), the lemma holds.

The following lemma can be derived by a simple induction.
\begin{Lemma}
For $J\neq0$, $\partial_{t}^{J}[G(v)]$ is a linear combination of terms
\begin{equation}
[D^{m}G](v)\partial_{t}^{\beta_{1}}v\partial_{t}^{\beta_{2}}v\cdots\partial_{t}^{\beta_{m}}v\quad\text{where}\quad 1\leq m\leq J,\quad \sum_{i=1}^{m}\beta_{m}=J.
\end{equation}
\end{Lemma}
Define

\begin{equation}
S_{1}(I,J):=\left\{\begin{array}{cc}
(\beta_{1},\cdots,\beta_{m},\alpha_{1},\cdots,\alpha_{m},\tilde{\alpha}_{1},\cdots,\tilde{\alpha}_{m},I_{1}
,I_{2},I_{3},J_{1},J_{2},J_{3}):\\
\max\{|I_{2}|+J_{2},|I_{3}|+J_{3}\}\leq|I|+J-1\\
\sum_{i=1}^{m}\beta_{i}+|\alpha_{i}|+|\tilde{\alpha}_{i}|+|I_{1}|+|I_{2}|+|I_{3}|+J_{1}+J_{2}+J_{3}\leq |I|+J
\end{array}
\right\},
\end{equation}
\begin{equation}
S_{2}(I,J):=\left\{\begin{array}{cc}
(\beta_{1},\cdots,\beta_{m},\alpha_{1},\cdots,\alpha_{m},\tilde{\alpha}_{1},\cdots,\tilde{\alpha}_{m},I_{1}
,I_{2},I_{3},J_{1},J_{2},J_{3}):\\
\sum_{i=1}^{m}\beta_{i}+|\alpha_{i}|+|\tilde{\alpha}_{i}|+|I_{1}|+|I_{2}|+|I_{3}|+J_{1}+J_{2}+J_{3}\leq |I|+J-1
\end{array}
\right\}
\end{equation}
and
\begin{equation}
S_{3}(I,J):=\left\{\begin{array}{cc}
(\beta_{1},\cdots,\beta_{m},\alpha_{1},\cdots,\alpha_{m},\tilde{\alpha}_{1},\cdots,\tilde{\alpha}_{m},I_{1}
,I_{2},I_{3},J_{1},J_{2},J_{3}):\\
\sum_{i=1}^{m}\beta_{i}+|\alpha_{i}|+|\tilde{\alpha}_{i}|+|I_{1}|+|I_{2}|+|I_{3}|+J_{1}+J_{2}+J_{3}\leq |I|+J-2
\end{array}
\right\}.
\end{equation}
Combing Lemmas 5.1-5.3, the following lemma holds obviously
\begin{Lemma}
For $J\geq 0$ and $|I|\geq0$, it holds that
\begin{equation}
\begin{split}
&\partial_{t}^{J}D^{I}\frac{Q(\varphi,e^{\alpha t}Q(\varphi,\varphi))}{2(1+e^{\alpha t}Q(\varphi,\varphi))}=\frac{Q(\varphi,e^{\alpha t}Q(\varphi,\partial_{t}^{J}D^{I}\varphi))}{(1+e^{\alpha t}
Q(\varphi,\varphi))}\\
&+\sum_{S_{1}(I,J)}\prod_{i=1}^{m}G(e^{\alpha t}Q(\varphi,\varphi))\partial_{t}^{\beta_{i}}(e^{\alpha t}Q(D^{\alpha_{i}}\varphi,D^{\tilde{\alpha}_{i}}\varphi))Q(\partial_{t}^{J_{1}}D^{I_{1}}\varphi,e^{\alpha t}Q(\partial_{t}^{J_{2}}D^{I_{2}}\varphi,\partial_{t}^{J_{3}}D^{I_{3}}\varphi))\\
&+\sum_{k=1}^{n}\sum_{S_{2}(I,J)}\prod_{i=1}^{m}G(e^{\alpha t}Q(\varphi,\varphi))\partial_{t}^{\beta_{i}}(e^{\alpha t}Q(D^{\alpha_{i}}\varphi,D^{\tilde{\alpha}_{i}}\varphi))Q(\partial_{t}^{J_{1}}D^{I_{1}}\varphi,e^{(\alpha-1)t}
\partial_{t}^{J_{2}}D^{I_{2}}\varphi_{k}\partial_{t}^{J_{3}}D^{I_{3}}\varphi_{k})\\
&+\sum_{k=1}^{n}\sum_{S_{2}(I,J)}\prod_{i=1}^{m}G(e^{\alpha t}Q(\varphi,\varphi))\partial_{t}^{\beta_{i}}(e^{\alpha t}Q(D^{\alpha_{i}}\varphi,D^{\tilde{\alpha}_{i}}\varphi))e^{(\alpha-1)t}(\partial_{t}^{J_{1}}D^{I_{1}}\varphi_{k})
\partial_{k}\left(Q(\partial_{t}^{J_{2}}D^{I_{2}}\varphi,\partial_{t}^{J_{3}}D^{I_{3}}\varphi)\right)\\
&+\sum_{k,j=1}^{n}\sum_{S_{3}(I,J)}\prod_{i=1}^{m}G(e^{\alpha t}Q(\varphi,\varphi))\partial_{t}^{\beta_{i}}(e^{\alpha t}Q(D^{\alpha_{i}}\varphi,D^{\tilde{\alpha_{i}}}\varphi))e^{(\alpha-2)t}(\partial_{t}^{J_{1}}D^{I_{1}}\varphi_{k})
(\partial_{t}^{J_{2}}D^{I_{2}}\varphi_{j})(\partial_{t}^{J_{3}}D^{I_{3}}\varphi_{j})\\
&:=\frac{Q(\varphi,e^{\alpha t}Q(\varphi,\partial_{t}^{J}D^{I}\varphi))}{(1+e^{\alpha t}
Q(\varphi,\varphi))}+R,
\end{split}
\end{equation}
where $G(e^{\alpha t}Q(\varphi,\varphi))$ denotes the set of smooth functions depending on $e^{\alpha t}Q(\varphi,\varphi)$ and $R$ stands for the remaining terms.
\end{Lemma}
\begin{Remark}
The term $\frac{Q(\varphi,e^{\alpha t}Q(\varphi,\varphi))}{1+e^{\alpha t}Q(\varphi,\varphi)}$ contains the highest order derivatives.
\end{Remark}
\begin{Remark}
From now on, without loss of generality, we assume that $\alpha_{0}+|\alpha|\leq \beta_{0}+|\beta|$ whenever they appear in $Q(\partial_{t}^{\alpha_{0}}D^{\alpha}\varphi, \partial_{t}^{\beta_{0}}D^{\beta}\varphi)$ simultaneously.
\end{Remark}
Now, we need to derive the energy inequality of the following equation
\begin{equation}
\begin{split}
\Box_{g}(\partial_{t}^{J}D^{I}\varphi)-(\partial_{t}^{J}D^{I}\varphi)_{t}&=e^{-t}(\sum_{i=1}^{n}\sum_{M=0}^{J-1}C_{M}\partial_{t}^{M}\partial_{i}^{2}D^{I}\varphi)\\
&+\frac{Q(\varphi,e^{\alpha t}Q(\varphi,\partial_{t}^{J}D^{I}\varphi))}{1+e^{\alpha t}Q(\varphi,\varphi)}+R,
\end{split}
\end{equation}
where $R$ is defined by (5.11).

The following lemma will play a key role in the proof of Theorem 1.2.
\begin{Lemma}
The following generalized energy inequality holds in the existence domain of the solution of Cauchy problem (1.9)-(1.10)
\begin{equation}
F(t)\leq 3\left(F(0)+\int_{0}^{t}(C (e^{-\frac{\tau}{2}}+e^{(\alpha-1)\tau}F(\tau))F(\tau)d\tau\right),
\end{equation}
provided the initial data is sufficiently small and $F(t)$ is defined by (4.7), $C$ is a constant depend only on $\alpha$ and $N$.
\end{Lemma}
{\it Proof.} We will prove the lemma by three steps:

Step 1: Energy estimates for (5.12).\\
Taking the vector field $V=-\partial_{t}$ and by Lemma 3.2, we have the following
\begin{equation}
E^{|I|,J}(t)-E^{|I|,J}(0)=\int_{0}^{t}\int_{\Sigma_{\tau}}(K^{V}(\partial_{\tau}^{J}D^{I}\varphi)+\Box_{g}(\partial_{\tau}^{J}D^{I}\varphi)V(\partial_{\tau}^{J}D^{I}\varphi))dVol_{\tau}d\tau.
\end{equation}
Which is equivalent to
\begin{equation}
\frac{d}{dt}E^{|I|,J}(t)=\int_{\Sigma_{t}}(K^{V}(\partial_{t}^{J}D^{I}\varphi(t))+\Box_{g}(\partial_{t}^{J}D^{I}\varphi)V(\partial_{t}^{J}D^{I}\varphi))dVol_{t}.
\end{equation}
Then, as (4.11), we get
\begin{equation}
\begin{split}
&\frac{d}{dt}f^{|I|,J}(t)\\
&=e^{\frac{2-n}{2}t}\frac{d}{dt}E^{|I|,J}(t)+\frac{2-n}{2}e^{\frac{2-n}{2}t}E^{|I|,J}(t)\\
&=e^{\frac{2-n}{2}t}\int_{\Sigma_{t}}\left(K^{V}(\partial_{t}^{J}D^{I}\varphi(t))+\Box_{g}(\partial_{t}^{J}D^{I}\varphi)V(\partial_{t}^{J}D^{I}\varphi)\right)dVol_{t}
+\frac{2-n}{2}e^{\frac{2-n}{2}t}E^{|I|,J}(t)\\
&=e^{\frac{2-n}{2}t}\int_{\Sigma_{t}}\left(e^{-t}(\sum_{i=1}^{n}\sum_{M=0}^{J-1}C_{M}\partial_{t}^{M}\partial_{i}^{2}D^{I}\varphi)+\frac{Q(\varphi,e^{\alpha t}Q(\varphi,\partial_{t}^{J}D^{I}\varphi))}{1+e^{\alpha t}Q(\varphi,\varphi)}+R\right)V(\partial_{t}^{J}D^{I}\varphi)dVol_{t}\\
&+e^{\frac{2-n}{2}t}\int_{\Sigma_{t}}\left(K^{V}(\partial_{t}^{J}D^{I}\varphi(t))+(\partial_{t}^{J}D^{I}\varphi)_{t}V(\partial_{t}^{J}D^{I}\varphi)\right)dVol_{t}
+\frac{2-n}{2}e^{\frac{2-n}{2}t}E^{|I|,J}(t)\\
&\leq e^{\frac{2-n}{2}t}\int_{\Sigma_{t}}\left(e^{-t}(\sum_{i=1}^{n}\sum_{M=0}^{J-1}C_{M}\partial_{t}^{M}\partial_{i}^{2}D^{I}\varphi)
+\frac{Q(\varphi,e^{\alpha t}Q(\varphi,\partial_{t}^{J}D^{I}\varphi))}{1+e^{\alpha t}Q(\varphi,\varphi)}+R\right)V(\partial_{t}^{J}D^{I}\varphi)dVol_{t}
\end{split}
\end{equation}
The last inequality holds according to Lemma 3.4. By (5.16), we have to estimate the following integral term containing the second order derivatives of $\partial_{t}^{J}D^{I}\varphi$.
\begin{equation}
\int_{\mathbb R^{n}}\frac{Q(\varphi,e^{\alpha t}Q(\varphi,\partial_{t}^{J}D^{I}\varphi))}{1+e^{\alpha t}Q(\varphi,\varphi)}V(\partial_{t}^{J}D^{I}\varphi)e^{\frac{n}{2}t}e^{\frac{2-n}{2}t}dx_{1}\cdots dx_{n}.
\end{equation}

Step 2: Estimates for (5.17).\\
By (1.6)
\begin{equation}
\begin{split}
&Q(\varphi,e^{\alpha t}Q(\varphi,\partial_{t}^{J}D^{I}\varphi))\\
&=-\varphi_{t}\partial_{t}\left[e^{\alpha t}\left(-\varphi_{t}(\partial_{t}^{J}D^{I}\varphi)_{t}+e^{-t}\sum_{i=1}^{n}\varphi_{i}
(\partial_{t}^{J}D^{I}\varphi)_{i}\right)\right]\\
&+e^{-t}\sum_{j=1}^{n}\varphi_{j}\partial_{j}\left[e^{\alpha t}\left(-\varphi_{t}(\partial_{t}^{J}D^{I}\varphi)_{t}+e^{-t}\sum_{i=1}^{n}\varphi_{i}
(\partial_{t}^{J}D^{I}\varphi)_{i}\right)\right]\\
&=e^{\alpha t}(\varphi_{t})^{2}(\partial_{t}^{J}D^{I}\varphi)_{tt}+e^{(\alpha-2)t}\sum_{i,j=1}^{n}\varphi_{i}\varphi_{j}(\partial_{t}^{J}D^{I}\varphi)_{ij}
-2e^{(\alpha-1)t}\sum_{i=1}^{n}\varphi_{t}\varphi_{i}(\partial_{t}^{J}D^{I}\varphi)_{ti}\\
&+\alpha e^{\alpha t}(\varphi_{t})^{2}(\partial_{t}^{J}D^{I}\varphi)_{t}+e^{\alpha t}\varphi_{t}\varphi_{tt}(\partial_{t}^{J}D^{I}\varphi)_{t}-\sum_{i=1}^{n}(\alpha-1)e^{(\alpha-1)t}\varphi_{t}\varphi_{i}(\partial_{t}^{J}D^{I}\varphi)_{i}\\
&-\sum_{i=1}^{n}e^{(\alpha-1)t}\varphi_{t}\varphi_{ti}
(\partial_{t}^{J}D^{I}\varphi)_{i}
-\sum_{j=1}^{n}e^{(\alpha-1)t}\varphi_{j}\varphi_{tj}(\partial_{t}^{J}D^{I}\varphi)_{t}
+e^{(\alpha-2)t}\sum_{i,j=1}^{n}\varphi_{j}\varphi_{ij}(\partial_{t}^{J}D^{I}\varphi)_{i}\\
&:=A+B+D+P,
\end{split}
\end{equation}
where
\begin{equation}
A=e^{\alpha t}(\varphi_{t})^{2}(\partial_{t}^{J}D^{I}\varphi)_{tt},
\end{equation}
\begin{equation}
B=e^{(\alpha-2)t}\sum_{i,j=1}^{n}\varphi_{i}\varphi_{j}(\partial_{t}^{J}D^{I}\varphi)_{ij}.
\end{equation}
\begin{equation}
D=-2e^{(\alpha-1)t}\sum_{i=1}^{n}\varphi_{t}\varphi_{i}(\partial_{t}^{J}D^{I}\varphi)_{it}
\end{equation}
and
\begin{equation}
\begin{split}
P&=\alpha e^{\alpha t}(\varphi_{t})^{2}(\partial_{t}^{J}D^{I}\varphi)_{t}+e^{\alpha t}\varphi_{t}\varphi_{tt}(\partial_{t}^{J}D^{I}\varphi)_{t}-\sum_{i=1}^{n}(\alpha-1)e^{(\alpha-1)t}\varphi_{t}\varphi_{i}(\partial_{t}^{J}D^{I}\varphi)_{i}\\
&-\sum_{i=1}^{n}e^{(\alpha-1)t}\varphi_{t}\varphi_{ti}
(\partial_{t}^{J}D^{I}\varphi)_{i}
-\sum_{j=1}^{n}e^{(\alpha-1)t}\varphi_{j}\varphi_{tj}(\partial_{t}^{J}D^{I}\varphi)_{t}+e^{(\alpha-2)t}\sum_{i,j=1}^{n}\varphi_{j}\varphi_{ij}(\partial_{t}^{J}D^{I}\varphi)_{i}\\
\end{split}
\end{equation}
Denote $\partial_{t}^{J}D^{I}\varphi$ by $v$, then by (5.17), (5.19) and integrating by parts
\begin{equation}
\begin{split}
&\int_{\mathbb R^{n}}\frac{A}{1+e^{\alpha t}Q(\varphi,\varphi)}(-v_{t})e^{\frac{nt}{2}}e^{\frac{2-n}{2}t}dx_{1}\cdots dx_{n}\\
&=\int_{\mathbb R^{n}}\frac{e^{\alpha t}(\varphi_{t})^{2}(v)_{tt}}{1+e^{\alpha t}Q(\varphi,\varphi)}(-v_{t})e^{\frac{nt}{2}}e^{\frac{2-n}{2}t}dx_{1}\cdots dx_{n}\\
&=-\frac{1}{2}\frac{d}{dt}\int_{\mathbb R^{n}}\frac{e^{\alpha t}\varphi_{t}^{2}}{1+e^{\alpha t}Q(\varphi,\varphi)}v_{t}^{2}e^{\frac{nt}{2}}e^{\frac{2-n}{2}t}dx_{1}\cdots dx_{n}+A_{1}+A_{2}+A_{3},
\end{split}
\end{equation}
where
\begin{equation}
A_{1}=\frac{1}{2}\int_{\mathbb R^{n}}
\frac{(\alpha+1)e^{\alpha t}\varphi_{t}^{2}}{1+e^{\alpha t}Q(\varphi,\varphi)}v_{t}^{2}e^{\frac{nt}{2}}e^{\frac{2-n}{2}t}dx_{1}\cdots dx_{n},
\end{equation}
\begin{equation}
A_{2}=\int_{\mathbb R^{n}}\frac{e^{\alpha t}\varphi_{t}\varphi_{tt}}{1+e^{\alpha t}Q(\varphi,\varphi)}v_{t}^{2}e^{\frac{nt}{2}}e^{\frac{2-n}{2}t}dx_{1}\cdots dx_{n}
\end{equation}
and
\begin{equation}
A_{3}=-\frac{1}{2}\int_{\mathbb R^{n}}\frac{e^{\alpha t}\varphi_{t}^{2}[e^{\alpha t}Q(\varphi,\varphi)]_{t}}{(1+e^{\alpha t}Q(\varphi,\varphi))^{2}}v_{t}^{2}e^{\frac{nt}{2}}e^{\frac{2-n}{2}t}dx_{1}\cdots dx_{n}.
\end{equation}
By (5.17), (5.20) and integrating by parts, it holds that
\begin{equation}
\begin{split}
&\int_{\mathbb R^{n}}\frac{B}{1+e^{\alpha t}Q(\varphi,\varphi)}(-v_{t})e^{\frac{nt}{2}}e^{\frac{2-n}{2}t}dx_{1}\cdots dx_{n}\\
&=\int_{\mathbb R^{n}}\frac{e^{(\alpha-2)t}\varphi_{i}\varphi_{j}v_{ij}}{1+e^{\alpha t}Q(\varphi,\varphi)}(-v_{t})e^{\frac{nt}{2}}e^{\frac{2-n}{2}t}dx_{1}\cdots dx_{n}\\
&=\frac{1}{2}\frac{d}{dt}\sum_{i,j=1}^{n}\int_{\mathbb R^{n}}\frac{e^{(\alpha-2)t}\varphi_{i}\varphi_{j}}{1+e^{\alpha t}Q(\varphi,\varphi)}v_{i}v_{j}e^{\frac{nt}{2}}e^{\frac{2-n}{2}t}dx_{1}\cdots dx_{n}+\sum_{\kappa=1}^{6}B_{\kappa},
\end{split}
\end{equation}
where
\begin{equation}
B_{1}=-\frac{1}{2}\sum_{i,j=1}^{n}\int_{\mathbb R^{n}}\frac{(\alpha-1)e^{(\alpha-2)t}\varphi_{i}\varphi_{j}}{1+e^{\alpha t}Q(\varphi,\varphi)}v_{i}v_{j}e^{\frac{nt}{2}}e^{\frac{2-n}{2}t}dx_{1}\cdots dx_{n},
\end{equation}
\begin{equation}
B_{2}=-\sum_{i,j=1}^{n}\int_{\mathbb R^{n}}\frac{e^{(\alpha-2)t}\varphi_{it}\varphi_{j}}{1+e^{\alpha t}Q(\varphi,\varphi)}v_{i}v_{j}e^{\frac{nt}{2}}e^{\frac{2-n}{2}t}dx_{1}\cdots dx_{n},
\end{equation}
\begin{equation}
B_{3}=\frac{1}{2}\sum_{i,j=1}^{n}\int_{\mathbb R^{n}}\frac{e^{(\alpha-2)t}\varphi_{i}\varphi_{j}[e^{\alpha t}Q(\varphi,\varphi)]_{t}}{(1+e^{\alpha t}Q(\varphi,\varphi))^{2}}v_{i}v_{j}e^{\frac{nt}{2}}e^{\frac{2-n}{2}t}dx_{1}\cdots dx_{n},
\end{equation}
\begin{equation}
B_{4}=\sum_{i,j=1}^{n}\int_{\mathbb R^{n}}\frac{e^{(\alpha-2)t}\varphi_{ij}\varphi_{j}}{1+e^{\alpha t}Q(\varphi,\varphi)}v_{i}v_{t}e^{\frac{nt}{2}}e^{\frac{2-n}{2}t}dx_{1}\cdots dx_{n},
\end{equation}
\begin{equation}
B_{5}=\sum_{i,j=1}^{n}\int_{\mathbb R^{n}}\frac{e^{(\alpha-2)t}\varphi_{i}\varphi_{jj}}{1+e^{\alpha t}Q(\varphi,\varphi)}v_{i}v_{t}e^{\frac{nt}{2}}e^{\frac{2-n}{2}t}dx_{1}\cdots dx_{n},
\end{equation}
and
\begin{equation}
B_{6}=-\sum_{i,j=1}^{n}\int_{\mathbb R^{n}}\frac{e^{(\alpha-2)t}\varphi_{i}\varphi_{j}[e^{\alpha t}Q(\varphi,\varphi)]_{j}}{(1+e^{\alpha t}Q(\varphi,\varphi))^{2}}v_{i}v_{t}e^{\frac{nt}{2}}e^{\frac{2-n}{2}t}dx_{1}\cdots dx_{n}.
\end{equation}
By (5.17), (5.21) and integrating by parts, we have
\begin{equation}
\begin{split}
&\int_{\mathbb R^{n}}\frac{D}{1+e^{\alpha t}Q(\varphi,\varphi)}(-v_{t})e^{\frac{nt}{2}}e^{\frac{2-n}{2}t}dx_{1}\cdots dx_{n}\\
&=\sum_{i=1}^{n}\int_{\mathbb R^{n}}\frac{2e^{(\alpha-1)t}\varphi_{t}\varphi_{i}}{1+e^{\alpha t}Q(\varphi,\varphi)}(v_{it}v_{t})e^{\frac{nt}{2}}e^{\frac{2-n}{2}t}dx_{1}\cdots dx_{n}\\
&:=D_{1}+D_{2}+D_{3},
\end{split}
\end{equation}
where
\begin{equation}
D_{1}=-\sum_{i=1}^{n}\int_{\mathbb R^{n}}\frac{e^{(\alpha-1)t}\varphi_{ti}\varphi_{i}}{1+e^{\alpha t}Q(\varphi,\varphi)}v_{t}^{2}e^{\frac{nt}{2}}e^{\frac{2-n}{2}t}dx_{1}\cdots dx_{n}.
\end{equation}
\begin{equation}
D_{2}=-\sum_{i=1}^{n}\int_{\mathbb R^{n}}\frac{e^{(\alpha-1)t}\varphi_{t}\varphi_{ii}}{1+e^{\alpha t}Q(\varphi,\varphi)}v_{t}^{2}e^{\frac{nt}{2}}e^{\frac{2-n}{2}t}dx_{1}\cdots dx_{n}
\end{equation}
and
\begin{equation}
D_{3}=\sum_{i=1}^{n}\int_{\mathbb R^{n}}\frac{e^{(\alpha-1)t}\varphi_{t}\varphi_{i}[e^{\alpha t}Q(\varphi,\varphi)]_{i}}{(1+e^{\alpha t}Q(\varphi,\varphi))^{2}}v_{t}^{2}e^{\frac{nt}{2}}e^{\frac{2-n}{2}t}dx_{1}\cdots dx_{n}
\end{equation}
Now, suppose that $N\geq n+6$ and $F(t)$ defined by (4.7) is small enough. According to Lemma 3.7 and Remark 3.2, we easily obtain for any $J+|I|\leq N$
\begin{equation}
|e^{\alpha t}Q(\varphi,\varphi)|=|e^{\alpha t}(-\varphi_{t}^{2}+\sum_{i=1}^{n}e^{-t}\varphi_{i}^{2})|\leq C e^{(\alpha-1)t}F(t)
\end{equation}
and
\begin{equation}
|\partial_{t}(e^{\alpha t}Q(\varphi,\varphi))|\leq |\alpha e^{\alpha t}Q(\varphi,\varphi)|
+|2e^{\alpha t}Q(\partial_{t}\varphi,\varphi)|+\sum_{i=1}^{n}|e^{(\alpha-1)t}\varphi_{i}^{2}|\leq C e^{(\alpha-1)t}F(t)
\end{equation}
then
\begin{equation}
\frac{1}{2}\leq1+e^{\alpha t}Q(\varphi,\varphi)\leq\frac{3}{2},
\end{equation}
\begin{equation}
|A_{1}|\leq \|\frac{(\alpha+1)e^{\alpha t}\varphi_{t}^{2}}{1+e^{\alpha t}Q(\varphi,\varphi)}\|_{L^{\infty}}F_{0}(t)
\leq Ce^{(\alpha-1)t}F(t)F_{0}(t).
\end{equation}
\begin{equation}
|A_{2}|\leq \|\frac{e^{\alpha t}\varphi_{t}\varphi_{tt}}{1+e^{\alpha t}Q(\varphi,\varphi)}\|_{L^{\infty}}F_{0}(t)
\leq Ce^{(\alpha-1)t}F(t)F_{0}(t).
\end{equation}
\begin{equation}
|A_{3}|\leq\|\frac{e^{\alpha t}\varphi_{t}^{2}[e^{\alpha t}Q(\varphi,\varphi)]_{t}}{(1+e^{\alpha t}Q(\varphi,\varphi))^{2}}\|_{L^{\infty}}F_{0}(t)\leq Ce^{2(\alpha-1)t}F^{2}(t)F_{0}(t).
\end{equation}
Similarly, we have
\begin{equation}
|B_{1}|\leq\|\frac{(\alpha-1)e^{(\alpha-1)t}\varphi_{i}\varphi_{j}}{1+e^{\alpha t}Q(\varphi,\varphi)}\|_{L^{\infty}}F_{1}(t)\leq Ce^{(\alpha-1)t}F(t)F_{1}(t).
\end{equation}
\begin{equation}
|B_{2}|\leq Ce^{(\alpha-\frac{3}{2})t}F(t)F_{1}(t),\quad |B_{3}|\leq C e^{2(\alpha-1)t}F^{2}(t)F_{1}(t).
\end{equation}
\begin{equation}
|B_{4}|\leq C e^{(\alpha-\frac{3}{2})t}F(t)F_{0}^{\frac{1}{2}}(t)F_{1}^{\frac{1}{2}}(t),\quad
|B_{5}|\leq C e^{(\alpha-\frac{3}{2})t}F(t)F_{0}^{\frac{1}{2}}(t)F_{1}^{\frac{1}{2}}(t).
\end{equation}
and
\begin{equation}
|B_{6}|\leq C e^{(2\alpha-\frac{5}{2})t}F^{2}(t)F_{0}^{\frac{1}{2}}(t)F_{1}^{\frac{1}{2}}(t).
\end{equation}
For $D_{i}$ $(i=1,2,3)$, we have
\begin{equation}
|D_{i}|\leq Ce^{(\alpha-\frac{3}{2})t}F(t)F_{0}(t)\quad i=1,2
\end{equation}
and
\begin{equation}
|D_{3}|\leq Ce^{(2\alpha-\frac{5}{2})t}F^{2}(t)F_{0}(t).
\end{equation}
At last, we estimate the term
\begin{equation}
\int_{\mathbb R^{n}}P(-v_{t})e^{\frac{nt}{2}}e^{\frac{2-n}{2}t}dx_{1}\cdots dx_{n}:=\sum_{\kappa=1}^{6}P_{\kappa},
\end{equation}
where $P_{\kappa}$ is defined orderly by the six parts of $P$. As before, we have
\begin{equation}
|P_{1}|\leq C e^{(\alpha-1)t}F(t)F_{0}(t),\quad |P_{2}|\leq C e^{(\alpha-1)t}F(t)F_{0}(t).
\end{equation}
\begin{equation}
|P_{3}|\leq C e^{(\alpha-1)t}F(t)F_{0}^{\frac{1}{2}}(t)F_{1}^{\frac{1}{2}}(t),\quad
|P_{4}|\leq C e^{(\alpha-\frac{3}{2})t}F(t)F_{0}^{\frac{1}{2}}(t)F_{1}^{\frac{1}{2}}(t)
\end{equation}
and
\begin{equation}
|P_{5}|\leq C e^{(\alpha-\frac{3}{2})t}F(t)F_{0}(t),\quad
|P_{6}|\leq C e^{(\alpha-1)t}F(t)F_{1}(t).
\end{equation}

In the last step, we estimate the remaining term
\begin{equation}
\int_{\mathbb R^{n}}R(-v_{t})e^{\frac{2-n}{2}t}e^{\frac{nt}{2}}dx_{1}\cdots dx_{n},
\end{equation}
where $R$ is defined by (5.11).

Step 3: Estimates for (5.54).\\
Before estimating (5.54), we expand every terms of $R$.

For $Q\left(\partial_{t}^{J_{1}}D^{I_{1}}\varphi,e^{\alpha t}Q(\partial_{t}^{J_{2}}D^{I_{2}}\varphi,\partial_{t}^{J_{3}}D^{I_{3}}\varphi)\right)$, it holds that
\begin{equation}
\begin{split}
&Q\left(\partial_{t}^{J_{1}}D^{I_{1}}\varphi,e^{\alpha t}Q(\partial_{t}^{J_{2}}D^{I_{2}}\varphi,\partial_{t}^{J_{3}}D^{I_{3}}\varphi)\right)\\
&=\alpha e^{\alpha t}(\partial_{t}^{J_{1}}D^{I_{1}}\varphi)_{t}(\partial_{t}^{J_{2}}D^{I_{2}}\varphi)_{t}(\partial_{t}^{J_{3}}D^{I_{3}}\varphi)_{t}+
e^{\alpha t}(\partial_{t}^{J_{1}}D^{I_{1}}\varphi)_{t}(\partial_{t}^{J_{2}}D^{I_{2}}\varphi)_{tt}(D^{I_{3}}\varphi)_{t}\\
&+e^{\alpha t}(\partial_{t}^{J_{1}}D^{I_{1}}\varphi)_{t}(\partial_{t}^{J_{2}}D^{I_{2}}\varphi)_{t}(\partial_{t}^{J_{3}}D^{I_{3}}\varphi)_{tt}
-\sum_{i=1}^{n}(\alpha-1)e^{(\alpha-1)t}(\partial_{t}^{J_{1}}D^{I_{1}}\varphi)_{t}(\partial_{t}^{J_{2}}D^{I_{2}}\varphi)_{i}(\partial_{t}^{J_{3}}D^{I_{3}}\varphi)_{i}\\
&-\sum_{i=1}^{n}e^{(\alpha-1)t}(\partial_{t}^{J_{1}}D^{I_{1}}\varphi)_{t}(\partial_{t}^{J_{2}}D^{I_{2}}\varphi)_{it}(\partial_{t}^{J_{3}}D^{I_{3}}\varphi)_{i}-
\sum_{i=1}^{n}e^{(\alpha-1)t}(\partial_{t}^{J_{1}}D^{I_{1}}\varphi)_{t}(\partial_{t}^{J_{2}}D^{I_{2}}\varphi)_{i}(\partial_{t}^{J_{3}}D^{I_{3}}\varphi)_{it}\\
&-\sum_{j=1}^{n}e^{(\alpha-1)t}(\partial_{t}^{J_{1}}D^{I_{1}}\varphi)_{j}(\partial_{t}^{J_{2}}D^{I_{2}}\varphi)_{tj}(\partial_{t}^{J_{3}}D^{I_{3}}\varphi)_{t}
-\sum_{j=1}^{n}e^{(\alpha-1)t}(\partial_{t}^{J_{1}}D^{I_{1}}\varphi)_{j}(\partial_{t}^{J_{2}}D^{I_{2}}\varphi)_{t}(\partial_{t}^{J_{3}}D^{I_{3}}\varphi)_{tj}\\
&+\sum_{i,j=1}^{n}e^{(\alpha-2)t}(\partial_{t}^{J_{1}}D^{I_{1}}\varphi)_{j}(\partial_{t}^{J_{2}}D^{I_{2}}\varphi)_{ij}(\partial_{t}^{J_{3}}D^{I_{3}}\varphi)_{i}+
\sum_{i,j=1}^{n}e^{(\alpha-2)t}(\partial_{t}^{J_{1}}D^{I_{1}}\varphi)_{j}(\partial_{t}^{J_{2}}D^{I_{2}}\varphi)_{i}(\partial_{t}^{J_{3}}D^{I_{3}}\varphi)_{ij}\\
&:=\sum_{\lambda=1}^{10}H_{\lambda},
\end{split}
\end{equation}
where $H_{\lambda}$ $(\lambda=1,\cdots,10)$ are defined orderly.

For $Q(\partial_{t}^{J_{1}}D^{I_{1}}\varphi,e^{(\alpha-1)t}\partial_{t}^{J_{2}}D^{I_{2}}\varphi_{i}
\partial_{t}^{J_{3}}D^{I_{3}}\varphi_{i})$,
 it holds that
 \begin{equation}
 \begin{split}
& Q(\partial_{t}^{J_{1}}D^{I_{1}}\varphi,e^{(\alpha-1)t}\partial_{t}^{J_{2}}
 D^{I_{2}}\varphi_{i}\partial_{t}^{J_{3}}D^{I_{3}}\varphi_{i})\\
 &=-(\alpha-1)e^{(\alpha-1)t}(\partial_{t}^{J_{1}}D^{I_{1}}\varphi)_{t}(\partial_{t}^{J_{2}}D^{I_{2}}\varphi_{i})
 (\partial_{t}^{J_{3}}D^{I_{3}}\varphi_{i})-e^{(\alpha-1)t}(\partial_{t}^{J_{1}}D^{I_{1}}\varphi)_{t}
 (\partial_{t}^{J_{2}}D^{I_{2}}\varphi_{i})_{t}(\partial_{t}^{J_{3}}D^{I_{3}}\varphi_{i})\\
 &-e^{(\alpha-1)t}(\partial_{t}^{J_{1}}D^{I_{1}}\varphi)_{t}
 (\partial_{t}^{J_{2}}D^{I_{2}}\varphi_{i})(\partial_{t}^{J_{3}}D^{I_{3}}\varphi_{i})_{t}+
 \sum_{j=1}^{n}e^{(\alpha-2)t}(\partial_{t}^{J_{1}}D^{I_{1}}\varphi)_{j}
 (\partial_{t}^{J_{2}}D^{I_{2}}\varphi_{i})_{j}(\partial_{t}^{J_{3}}D^{I_{3}}\varphi_{i})\\
 &+\sum_{j=1}^{n}e^{(\alpha-2)t}(\partial_{t}^{J_{1}}D^{I_{1}}\varphi)_{j}
 (\partial_{t}^{J_{2}}D^{I_{2}}\varphi_{i})(\partial_{t}^{J_{3}}D^{I_{3}}\varphi_{i})_{j}\\
 &:=\sum_{\lambda=1}^{5}O_{\lambda}.
 \end{split}
 \end{equation}

 For $e^{(\alpha-1)t}(\partial_{t}^{J_{1}}D^{I_{1}}\varphi_{i})
 \partial_{i}\left(Q(\partial_{t}^{J_{2}}D^{I_{2}}\varphi,\partial_{t}^{J_{3}}D^{I_{3}}\varphi)\right)$, it holds that
 \begin{equation}
 \begin{split}
 &e^{(\alpha-1)t}(\partial_{t}^{J_{1}}D^{I_{1}}\varphi_{i})
 \partial_{i}\left(Q(\partial_{t}^{J_{2}}D^{I_{2}}\varphi,\partial_{t}^{J_{3}}D^{I_{3}}\varphi)\right)\\
 &=-e^{(\alpha-1)t}(\partial_{t}^{J_{1}}D^{I_{1}}\varphi_{i})(\partial_{t}^{J_{2}}D^{I_{2}}\varphi_{i})_{t}
 (\partial_{t}^{J_{3}}D^{I_{3}}\varphi)_{t}-e^{(\alpha-1)t}(\partial_{t}^{J_{1}}D^{I_{1}}\varphi_{i})
 (\partial_{t}^{J_{2}}D^{I_{2}}\varphi)_{t}
 (\partial_{t}^{J_{3}}D^{I_{3}}\varphi_{i})_{t}\\
 &+\sum_{j=1}^{n}e^{(\alpha-2)t}(\partial_{t}^{J_{1}}D^{I_{1}}\varphi_{i})
 (\partial_{t}^{J_{2}}D^{I_{2}}\varphi_{i})_{j}(\partial_{t}^{J_{3}}D^{I_{3}}\varphi)_{j}
 +\sum_{j=1}^{n}e^{(\alpha-2)t}(\partial_{t}^{J_{1}}D^{I_{1}}\varphi_{i})
 (\partial_{t}^{J_{2}}D^{I_{2}}\varphi)_{j}(\partial_{t}^{J_{3}}D^{I_{3}}\varphi_{i})_{j}\\
 &:=\sum_{\lambda=1}^{4}Q_{\lambda}.
 \end{split}
 \end{equation}

 Denote
 \begin{equation}
 X=e^{(\alpha-2)t}(\partial_{t}^{J_{1}}
 D^{I_{1}}\varphi_{i})(\partial_{t}^{J_{2}}D^{I_{2}}\varphi_{j})(\partial_{t}^{J_{3}}D^{I_{3}}\varphi_{j}).
 \end{equation}

 At last, for $\partial_{t}^{\beta_{i}}(e^{\alpha t}Q(D^{\alpha_{i}}\varphi,D^{\tilde{\alpha}_{i}}\varphi))$, it holds that
 \begin{equation}
 \begin{split}
 &\partial_{t}^{\beta_{i}}(e^{\alpha t}Q(D^{\alpha_{i}}\varphi,D^{\tilde{\alpha}_{i}}\varphi))\\
 &=e^{\alpha t}\sum_{\beta_{i1}+\beta_{i2}\leq\beta_{i}}\partial_{t}^{\beta_{i1}}(D^{\alpha_{i}}\varphi)_{t}\partial_{t}^{\beta_{i2}}(D^{\tilde{\alpha_{i}}}\varphi)_{t}
 +e^{(\alpha-1)t}\sum_{j=1}^{n}\sum_{\beta_{i1}+\beta_{i2}\leq\beta_{i}}\partial_{t}^{\beta_{i1}}
 (D^{\alpha_{i}}\varphi)_{j}\partial_{t}^{\beta_{i2}}(D^{\tilde{\alpha}_{i}}\varphi)_{j}\\
 &:=\sum_{\lambda=1}^{2}Y_{\lambda},
 \end{split}
 \end{equation}
 where we have omit the constant coefficients, which do not affect the main result.
We estimate (5.54) by the following four cases according to the index.\\

Case I: when $\sum_{i=1}^{m}\beta_{i}+\alpha_{i}+\tilde{\alpha_{i}}=0$, we have
\begin{equation}
\begin{split}
&\int_{\mathbb R^{n}}R(-v_{t})e^{\frac{2-n}{2}t}e^{\frac{nt}{2}}dx_{1}\cdots dx_{n}\\
&\leq\int_{\mathbb R^{n}}|G\left(e^{\alpha t}Q(\varphi,\varphi)\right)|(\sum_{\lambda=1}^{10}\sum_{S_{1}(I,J)}|H_{\lambda}|+\sum_{\lambda=1}^{5}\sum_{S_{2}(I,J)}|O_{\lambda}|)
e^{\frac{2-n}{2}t}e^{\frac{nt}{2}}dx_{1}\cdots dx_{n}\\
&+\int_{\mathbb R^{n}}|G\left(e^{\alpha t}Q(\varphi,\varphi)\right)(\sum_{\lambda=1}^{4}\sum_{S_{2}(I,J)}|Q_{\lambda}|+\sum_{S_{3}(I,J)}|X|)e^{\frac{2-n}{2}t}e^{\frac{nt}{2}}dx_{1}\cdots dx_{n}
\end{split}
\end{equation}
 We will use the following principle to estimate these product terms here and hereafter.

{\bf Principle}: Since $\max\{J_{2}+|I_{2}|,J_{3}+|I_{3}|\}\leq |I|+J-1$ and $\sum_{i=1}^{3}(|I_{i}|+J_{i})\leq |I|+J$, there must be at most one term that exceeds $\frac{|I|}{2}$, we use $L^{2}$ norm to control this term and use $L^{\infty}$ norm to control other terms.

According the above principle and step 2, when $|I_{1}|+J_{1}\geq\frac{|I|+J}{2}$, by Lemma 3.9, since $\frac{N}{2}+2+\frac{n}{2}+1\leq N$, it holds that
\begin{equation}
\int_{\mathbb R^{n}}|G(e^{\alpha t}Q(\varphi,\varphi))H_{i}v_{t}|e^{\frac{nt}{2}}e^{\frac{2-n}{2}t}dx_{1}\cdots dx_{n}
\leq Ce^{(\alpha-1)t}F(t)F_{0}(t)\quad i=1,\cdots,4
\end{equation}
\begin{equation}
\int_{\mathbb R^{n}}|G(e^{\alpha t}Q(\varphi,\varphi))H_{i}v_{t}|e^{\frac{nt}{2}}e^{\frac{2-n}{2}t}dx_{1}\cdots dx_{n}
\leq Ce^{(\alpha-\frac{3}{2})t}F(t)F_{0}(t)\quad i=5,6
\end{equation}
\begin{equation}
\int_{\mathbb R^{n}}|G(e^{\alpha t}Q(\varphi,\varphi))H_{i}v_{t}|e^{\frac{nt}{2}}e^{\frac{2-n}{2}t}dx_{1}\cdots dx_{n}
\leq Ce^{(\alpha-\frac{3}{2})t}F(t)F_{0}^{\frac{1}{2}}(t)F_{1}^{\frac{1}{2}}(t)\quad i=7,\cdots,10
\end{equation}
\begin{equation}
\int_{\mathbb R^{n}}|G(e^{\alpha t}Q(\varphi,\varphi))O_{i}v_{t}|e^{\frac{nt}{2}}e^{\frac{2-n}{2}t}dx_{1}\cdots dx_{n}
\leq Ce^{(\alpha-1)t}F(t)^{2}\quad i=1
\end{equation}
\begin{equation}
\int_{\mathbb R^{n}}|G(e^{\alpha t}Q(\varphi,\varphi))O_{i}v_{t}|e^{\frac{nt}{2}}e^{\frac{2-n}{2}t}dx_{1}\cdots dx_{n}
\leq Ce^{(\alpha-\frac{3}{2})t}F(t)^{2}\quad i=2,3,4,5
\end{equation}
\begin{equation}
\int_{\mathbb R^{n}}|G(e^{\alpha t}Q(\varphi,\varphi))Q_{i}v_{t}|e^{\frac{nt}{2}}e^{\frac{2-n}{2}t}dx_{1}\cdots dx_{n}
\leq Ce^{(\alpha-\frac{3}{2})t}F(t)^{2}\quad i=1,2,3,4
\end{equation}
and
\begin{equation}
\int_{\mathbb R^{n}}|G(e^{\alpha t}Q(\varphi,\varphi))X v_{t}|e^{\frac{nt}{2}}e^{\frac{2-n}{2}t}dx_{1}\cdots dx_{n}
\leq Ce^{(\alpha-\frac{3}{2})t}F(t)^{2}.
\end{equation}

When $|I_{3}|+J_{3}\geq\frac{|I|+J}{2}$, we have
\begin{equation}
\int_{\mathbb R^{n}}|G(e^{\alpha t}Q(\varphi,\varphi))H_{i}v_{t}|e^{\frac{nt}{2}}e^{\frac{2-n}{2}t}dx_{1}\cdots dx_{n}
\leq Ce^{(\alpha-1)t}F(t)F_{0}(t)\quad i=1,\cdots,3
\end{equation}
\begin{equation}
\int_{\mathbb R^{n}}|G(e^{\alpha t}Q(\varphi,\varphi))H_{i}v_{t}|e^{\frac{nt}{2}}e^{\frac{2-n}{2}t}dx_{1}\cdots dx_{n}
\leq Ce^{(\alpha-1)t}F(t)F_{0}^{\frac{1}{2}}(t)F_{1}^{\frac{1}{2}}(t)\quad i=4
\end{equation}
\begin{equation}
\int_{\mathbb R^{n}}|G(e^{\alpha t}Q(\varphi,\varphi))H_{i}v_{t}|e^{\frac{nt}{2}}e^{\frac{2-n}{2}t}dx_{1}\cdots dx_{n}
\leq Ce^{(\alpha-\frac{3}{2})t}F(t)F_{0}^{\frac{1}{2}}(t)F_{1}^{\frac{1}{2}}(t)\quad i=5,9,10
\end{equation}
\begin{equation}
\int_{\mathbb R^{n}}|G(e^{\alpha t}Q(\varphi,\varphi))H_{i}v_{t}|e^{\frac{nt}{2}}e^{\frac{2-n}{2}t}dx_{1}\cdots dx_{n}
\leq Ce^{(\alpha-\frac{3}{2})t}F(t)F_{0}(t)\quad i=6,7,8
\end{equation}
\begin{equation}
\int_{\mathbb R^{n}}|G(e^{\alpha t}Q(\varphi,\varphi))O_{i}v_{t}|e^{\frac{nt}{2}}e^{\frac{2-n}{2}t}dx_{1}\cdots dx_{n}
\leq Ce^{(\alpha-1)t}F(t)^{2}\quad i=1
\end{equation}
\begin{equation}
\int_{\mathbb R^{n}}|G(e^{\alpha t}Q(\varphi,\varphi))O_{i}v_{t}|e^{\frac{nt}{2}}e^{\frac{2-n}{2}t}dx_{1}\cdots dx_{n}
\leq Ce^{(\alpha-\frac{3}{2})t}F(t)^{2}\quad i=2,3,4,5
\end{equation}
\begin{equation}
\int_{\mathbb R^{n}}|G(e^{\alpha t}Q(\varphi,\varphi))Q_{i}v_{t}|e^{\frac{nt}{2}}e^{\frac{2-n}{2}t}dx_{1}\cdots dx_{n}
\leq Ce^{(\alpha-\frac{3}{2})t}F(t)^{2}\quad i=1,2,3,4
\end{equation}
and
\begin{equation}
\int_{\mathbb R^{n}}|G(e^{\alpha t}Q(\varphi,\varphi))Xv_{t}|e^{\frac{nt}{2}}e^{\frac{2-n}{2}t}dx_{1}\cdots dx_{n}
\leq Ce^{(\alpha-\frac{3}{2})t}F(t)^{2}
\end{equation}

In the above calculations, we have used the fact that $G(e^{\alpha t}Q(\varphi,\varphi))$ is uniformly bounded, provided  $F(t)$ is sufficiently small and the uniform constant $C$ depends only on $\alpha$ and $N$. Thus, by (5.61)-(5.75), we conclude that in this case
\begin{equation}
\int_{\mathbb R^{n}}R(-v_{t})e^{\frac{2-n}{2}t}e^{\frac{nt}{2}}dx_{1}\cdots dx_{n}
\leq C e^{(\alpha-1)t}F^{2}(t).
\end{equation}

Case II: when $\sum_{i=1}^{m}\beta_{i}+\alpha_{i}+\tilde{\alpha}_{i}=j>0$, and $max\{|I_{1}|+J_{1},|I_{2}|+J_{2},|I_{3}|+J_{3}\}\geq\frac{|I|+J}{2}$.\\
In this case, we have to do some additional estimates on the $Y_{\lambda}$ terms by the $L^{\infty}$ norm and get some additional decay. By direct calculations, we have
\begin{equation}
|Y_{\lambda}|\leq C e^{(\alpha-1)t}F(t).
\end{equation}
Then, by (5.76) and (5.77), we conclude that
\begin{equation}
\int_{\mathbb R^{n}}R(-v_{t})e^{\frac{2-n}{2}t}e^{\frac{nt}{2}}dx_{1}\cdots dx_{n}
\leq C e^{2(\alpha-1)t}F^{3}(t).
\end{equation}

Case III: when $\sum_{i=1}^{m}\beta_{i}+\alpha_{i}+\tilde{\alpha_{i}}=j>0$, there exists some $i$ such that $|\alpha_{i}|+|\tilde{\alpha_{i}}|+\beta_{i}\geq\frac{|I|+J}{2}$.\\
In this case, we use the principle of Case I, we estimate (5.55)-(5.58) by $L^{\infty}$ norm, by direct calculations, we have
\begin{equation}
|H_{i}|\leq C e^{(\alpha-\frac{3}{2})t}F^{\frac{3}{2}}(t),\qquad i=1,\cdots,4
\end{equation}
\begin{equation}
|H_{i}|\leq C e^{(\alpha-2)t}F^{\frac{3}{2}}(t),\qquad i=5,\cdots,10
\end{equation}
\begin{equation}
|O_{i}|\leq C e^{(\alpha-\frac{3}{2})t}F^{\frac{3}{2}}(t),\qquad i=1
\end{equation}
\begin{equation}
|O_{i}|\leq C e^{(\alpha-2)t}F^{\frac{3}{2}}(t),\qquad i=2,3,4,5
\end{equation}
\begin{equation}
|Q_{i}|\leq C e^{(\alpha-2)t}F^{\frac{3}{2}}(t),\qquad i=1,2,3,4
\end{equation}
and
\begin{equation}
|X|\leq C e^{(\alpha-2)t}F^{\frac{3}{2}}(t).
\end{equation}
For the term $\partial_{t}^{\beta_{i}}(e^{\alpha t}Q(D^{\alpha_{i}}\varphi,D^{\tilde{\alpha_{i}}}\varphi))$, where $|\alpha_{i}|+\beta_{i}+|\tilde{\alpha_{i}}|\geq\frac{|I|+J}{2}$, we have
\begin{equation}
\int_{\mathbb R^{n}}|Y_{i}v_{t}|e^{\frac{nt}{2}}e^{\frac{2-n}{2}t}dx_{1}\cdots dx_{n}
\leq C e^{(\alpha-\frac{1}{2})t}F^{\frac{1}{2}}(t)F(t).
\end{equation}
Then, by (5.79)-(5.85), we conclude that in this case, we have
\begin{equation}
\int_{\mathbb R^{n}}R(-v_{t})e^{\frac{2-n}{2}t}e^{\frac{nt}{2}}dx_{1}\cdots dx_{n}
\leq C e^{2(\alpha-1)t}F^{3}(t).
\end{equation}
Case IV: $max\{|I_{1}|+J_{1},|I_{2}|+J_{2},|I_{3}|+J_{3},\beta_{i}+|\alpha_{i}|+|\tilde{\alpha_{i}}|\}\leq\frac{|I|+J}{2}$.\\
This case is easy to deal with, we can estimate it by Case II and Case III.

Combing Cases I-IV, it holds that
\begin{equation}
\int_{\mathbb R^{n}}R(-v_{t})e^{\frac{2-n}{2}t}e^{\frac{nt}{2}}dx_{1}\cdots dx_{n}
\leq C e^{(\alpha-1)t}F^{2}(t)
\end{equation}
 provided $F(t)$ is small enough.

By the above three steps and (4.13), summing all $|I|$, $J$ satisfying $|I|+J\leq N$, we have
\begin{equation}
\begin{split}
&\frac{d}{dt}(F(t)+\frac{1}{2}\int_{\mathbb R^{n}}\frac{e^{\alpha t}\varphi_{t}^{2}}{1+e^{\alpha t}Q(\varphi,\varphi)}v_{t}^{2}e^{\frac{nt}{2}}e^{\frac{2-n}{2}t}dx_{1}\cdots dx_{n}\\
&-\frac{1}{2}\sum_{i,j=1}^{n}\int_{\mathbb R^{n}}\frac{e^{(\alpha-2)t}\varphi_{i}\varphi_{j}}{1+e^{\alpha t}Q(\varphi,\varphi)}v_{i}v_{j}e^{\frac{nt}{2}}e^{\frac{2-n}{2}t}dx_{1}\cdots dx_{n})\\
&\leq C_{N}(e^{-\frac{t}{2}}+e^{(\alpha-1)t}F(t))F(t).
\end{split}
\end{equation}
The following holds by direct calculations
\begin{equation}
|\frac{1}{2}\int_{\mathbb R^{n}}\frac{e^{\alpha t}\varphi_{t}^{2}}{1+e^{\alpha t}Q(\varphi,\varphi)}v_{t}^{2}e^{\frac{nt}{2}}e^{\frac{2-n}{2}t}dx_{1}\cdots dx_{n}|
\leq Ce^{(\alpha-1)t}F^{2}(t)\leq\frac{1}{4}F(t)
\end{equation}
and
\begin{equation}
|\frac{1}{2}\sum_{i,j=1}^{n}\int_{\mathbb R^{n}}\frac{e^{(\alpha-2)t}\varphi_{i}\varphi_{j}}{1+e^{\alpha t}Q(\varphi,\varphi)}v_{i}v_{j}e^{\frac{nt}{2}}e^{\frac{2-n}{2}t}dx_{1}\cdots dx_{n})|
\leq Ce^{(\alpha-1)t}F^{2}(t)\leq\frac{1}{4}F(t),
\end{equation}
provided $F(t)$ is sufficiently small. By (5.88)-(5.96), we obtain
$$
F(t)\leq 3\left(F(0)+\int_{0}^{t}(C_{N} (e^{-\frac{\tau}{2}}+e^{(\alpha-1)\tau}F(\tau))F(\tau)d\tau\right).
$$
Thus, the lemma is proved.

{\it Proof of Theorem 1.2.}\\
As in the last section, what we have to do is to close the boot-strap assumption in the existence domain $[0,T)$. Define $E(t)$ as Section 4, we will prove that for any $t\in[0,T)$, $F(t)\leq 2A\epsilon$ implies $F(t)\leq A\epsilon$, provided $\epsilon$ is small enough.
By (4.4), Lemma 5.2 and Gronwall's lemma, we have
\begin{equation}
F(t)\leq 3F(0)e^{\int_{0}^{t}C(e^{-\frac{\tau}{2}}+ e^{(\alpha-1)\tau}(2A\epsilon))d\tau}\leq\frac{3C_{0}\epsilon}{2}e^{\int_{0}^{t}C_{N}(e^{-\frac{\tau}{2}}+ e^{(\alpha-1)\tau}(2A\epsilon))d\tau}.
\end{equation}
By (5.91), we have to choose $\epsilon_{0}$ sufficiently small such that
\begin{equation}
e^{\int_{0}^{t}C_{N}(e^{-\frac{\tau}{2}}+ e^{(\alpha-1)\tau}(2A\epsilon))d\tau}\leq\frac{2A}{3C_{0}}.
\end{equation}
When $\alpha<1$, we have
\begin{equation}
e^{\int_{0}^{t}C_{N}(e^{-\frac{\tau}{2}}+ e^{(\alpha-1)\tau}(2A\epsilon))d\tau}\leq e^{\int_{0}^{\infty}C_{N}(e^{-\frac{\tau}{2}}+ e^{(\alpha-1)\tau}(2A\epsilon))d\tau}\leq\frac{2A}{3C_{0}}
\end{equation}
Provided $\epsilon\leq\frac{1-\alpha}{2A}$ and $A\geq \frac{3C_{0}e^{3C_{N}}}{2}$ is sufficiently large.

When $\alpha=1$, by (5.92), we have
\begin{equation}
T(\epsilon)=\frac{\ln(C)}{2AC_{N}\epsilon},
\end{equation}
where $C=\frac{2A}{3C_{0}e^{2C_{N}}}$ is a positive constant.
Combining (5.93) and (5.94), Theorem 1.2 holds.
\begin{Remark}
From the procedure we used to prove the main Theorem 1.2, we see that the role played by the term $e^{\alpha t}$ is important on the lifespan of the solution in curved spacetime.
\end{Remark}
\section{Discussions}
In this paper, we use the method of vector fields to prove the well-posedness of nonlinear wave equations in a curved spacetime. By this method, we get the exponential decay for the spacetime derivatives, which is important to the nonlinear problems. Since de Sitter spacetime is a special case of the Robertson-Walker spacetime, whose metric has the following form
\begin{equation}
ds^{2}=-dt^{2}+a^{2}(t)(\sum_{i=1}^{n}dx^{2}_{i}),
\end{equation}
where $a(t)$ is an appropriate factor, we can get the similar results for the general Robertson-Walker spacetime depending on the choice of $a(t)$.

Since de Sitter spacetime is a special Lorentz spacetime, which is conformally flat. There maybe some more decay properties should be explored to improve the results of the present paper, such as conformal inequality, Morawetz inequality and so on, which have played a vital role in the Minkowski spacetime $\mathbb R^{1+n}$. We will discuss these inequalities in our forthcoming paper.
%

At last, it is interesting to study the wave equation satisfying null condition in the domain containing the black hole region on Schwarzschild-de Sitter spacetime and this maybe give a clue to the study of the stability of this spacetime as solutions to the Einstein equations, which is still an open problem in gravitational physics.
%
%

\vskip 6mm

\noindent{\Large {\bf Acknowledgements.}} This work was supported in part by the NNSF of China (Grant Nos.: 11271323,91330105) and the Zhejiang Provincial Natural Science Foundation of China (Grant No.: LZ13A010002).

\end{document}